
\input epsf.tex
\input amssym.def
\input amssym
\magnification=1000
\baselineskip = 0.21truein
\lineskiplimit = 0.01truein
\lineskip = 0.01truein
\vsize = 8.7truein
\voffset = 0.2truein
\parskip = 0.10truein
\parindent = 0.3truein
\settabs 12 \columns
\hsize = 6.0truein
\hoffset = 0.1truein

\setbox\strutbox=\hbox{%
\vrule height .708\baselineskip
depth .292\baselineskip
width 0pt}
\font\caps=cmcsc10
\font\bigtenrm=cmr10 at 14pt

\def\sqr#1#2{{\vcenter{\vbox{\hrule height.#2pt
\hbox{\vrule width.#2pt height#1pt \kern#1pt
\vrule width.#2pt}
\hrule height.#2pt}}}}
\def\square{\mathchoice\sqr46\sqr46\sqr{3.1}6\sqr{2.3}4}

\centerline{\bigtenrm LINKS WITH SPLITTING NUMBER ONE}
\tenrm
\vskip 14pt
\centerline{MARC LACKENBY\footnote{\dag}{Partially supported by EPSRC grant EP/R005125/1}}
\vskip 12pt
\centerline{Mathematical Institute, University of Oxford,}
\centerline{Radcliffe Observatory Quarter, Woodstock Road, Oxford OX2 6GG, United Kingdom.}
\centerline{Email: lackenby@maths.ox.ac.uk \ \ \ \ Tel: 01865 273561}
\vskip 18pt

\centerline{\caps Abstract}
\vskip 6pt
\noindent
{\narrower
We provide an algorithm to determine whether a link $L$ admits a crossing change that turns it into a split link, under some fairly mild hypotheses on $L$. The algorithm also provides a complete list of all such crossing changes. It can therefore also determine whether the unlinking number of $L$ is $1$. \par}

\vskip 18pt
\centerline{\caps 1. Introduction}
\vskip 6pt

One of the knot invariants that is least well understood is unknotting number. This is defined to be
the minimal number of crossing changes that one can apply to some diagram of the knot
in order to unknot it. For any given diagram of a knot $K$, it is of course easy to determine
the minimal number of crossing changes that one can apply to it in order to unknot it, by using one of the several known
algorithms to detect the unknot. However, one has no guarantee in general that there is not
some more complicated diagram of $K$ that can be unknotted using fewer crossing changes.
Many techniques have been developed to find lower bounds on the unknotting number
of a knot, for example, using the Alexander module [37], the Goeritz form [25, 47], 
gauge theory [48] and Heegaard Floer homology [31, 38, 39]. However, no known technique
is perfect, and in fact there are many explicit knots for which the unknotting number is
not known [26]. A satisfactory resolution will only be found when an algorithm that determines 
the unknotting number of a knot is discovered. But this appears
to be a very long way off. In fact, it is conceivable that no such algorithm exists.
It is not even known whether one can decide algorithmically
whether a knot has unknotting number one. 

In this paper, we explore some natural generalisations of unknotting number to links
with more than one component. One might consider the {\sl unlinking number}
$u(L)$ of a link L, which is the minimal number of crossing changes required
to turn it into the unlink. But it turns out that it is just as natural to consider the
{\sl splitting number} $s(L)$, which is the minimal number of crossing changes
required to turn it into a split link. (A link is {\sl split} if there is an embedded 2-sphere disjoint
from the link with link components on both sides.) Some authors [6] have also analysed a variant
of splitting number, where one only considers crossing changes between distinct
components of the link. The minimal number of such crossing changes required
to create a split link we denote by $s_d(L)$ (where $d$ stands for `distinct'). Other
authors [27] have also required the resulting link to be {\sl totally split}, which means
that there is a union of disjoint balls containing the link, such that each ball contains
a single component of the link in its interior. We say that the {\sl total splitting number}
$ts(L)$ is the minimal number of crossing changes required to make the link
totally split. Again, one can consider only crossing changes between distinct link components
and we denote the resulting variant of total splitting number by $ts_d(L)$.
Our main result is that, under some fairly mild
hypotheses, there are algorithms to determine whether any of these quantities
is 1 for a given link.

\noindent {\bf Theorem 1.1.} {\sl There is an algorithm to solve the following problem.
The input to the algorithm is a link $L$ in $S^3$, given either by a diagram or by a triangulation of 
$S^3$ with $L$ as a specified subcomplex.
The link $L$ must be hyperbolic and 2-string prime. It is required to have at least two components
and if it has exactly two components, these must have zero linking number.
The output is an answer to each of the following questions:
\vskip -6pt
\item{(i)} Is $u(L) = 1$?
\vskip -6pt
\item{(ii)} Is $s(L) = 1$?
\vskip -6pt
\item{(iii)} Is $s_d(L) = 1$?
\vskip -6pt
\item{(iv)} Is $ts(L) = 1$?
\vskip -6pt
\item{(v)} Is $ts_d(L) = 1$?

}

Recall that a link $L$ is {\sl 2-string prime} if, for each 2-sphere $S$ in $S^3$ that intersects $L$ transversely in four points,
$S - L$ admits a compression disc in the complement of $L$. When $L$ is hyperbolic, then it is 2-string
prime if and only if its branched double cover is hyperbolic or a small Seifert fibre space. (This is explained
in Section 4.) This condition can be readily verified both in theory [22, 46] and in practice [52].

The linking number hypothesis when $L$ has two components is a slightly unfortunate one.
However, it is not as restrictive as it first may seem. When $s(L) = 1$, then
the two components of $L$ must have linking number $0$, $1$ or $-1$, for the following reason.
When both components of $L$ are involved in the crossing change, then the
linking number is $\pm 1$, since the crossing change alters the linking number
by one. On the other hand, when the crossing change moves some
component of $L$ through itself, then this does not change the linking number,
and so this must be zero.

The most notable hypothesis in Theorem 1.1 is that the link $L$ has more than one
component. As mentioned above, it remains the case that there is no known algorithm
to decide whether a knot has unknotting number one.

Whenever one has an algorithmic result such as the one presented above, a finiteness theorem tends to
come for free. In this case, the question that we can address is: are there
only finitely many ways to split a link by a crossing change? Of course, one needs
a way to compare two crossing changes, which may occur in different diagrams.
There is a natural method of doing this using surgery. Given any crossing in some diagram of $L$,
we may encircle the two sub-arcs of $L$ near the crossing by a simple closed curve $C$,
as shown in Figure 1. This bounds an embedded disc $D$ such that $D \cap L$ is two points
in the interior of $D$. Such a disc $D$ is called a {\sl crossing disc}. The boundary curve
$C$ of a crossing disc is called a {\sl crossing circle}. In the interior of the crossing disc $D$, there is an embedded arc joining
the two points of $D \cap L$. This is the {\sl crossing arc} associated with $D$. Changing the crossing is
achieved by $\pm 1$ surgery along $C$. We say that two crossing changes are 
{\sl equivalent} if their associated crossing circles are ambient isotopic in the
complement of $L$ and the associated surgery coefficients are equal.

\noindent {\bf Theorem 1.2.} {\sl Let $L$ be as in Theorem 1.1. If $s(L) = 1$, then,
up to equivalence, there are only finitely many ways to turn $L$ into a split link by 
performing a crossing change. In fact, if $t$ is
the number of tetrahedra in a triangulation of $S^3$ with $L$ as a subcomplex, then
the number of distinct ways of creating a split link from $L$ by a crossing
change is at most $k^t$, for some universal computable constant $k$.
Hence, the number of ways is at most $k^{24 c(L)}$ where $c(L)$ is the
crossing number of $L$. Moreover, there is an algorithm to find all these
crossing changes.}

\vskip 12pt
\centerline{
\epsfxsize = 3in
\epsfbox{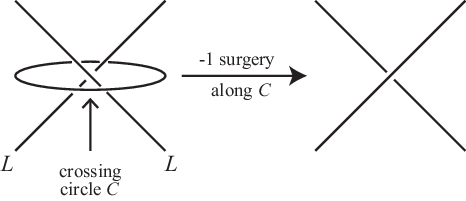}
}
\vskip 6pt
\centerline{Figure 1}

Note that a version of Theorem 1.2 also holds for the other variants of splitting number and unlinking number
discussed above. Indeed if any of $u(L)$, $s_d(L)$, $ts(L)$ or $ts_d(L)$ is equal to one, then necessarily $s(L) = 1$,
provided $L$ is non-split. So, in any of these cases, Theorem 1.2 also provides a finiteness result on the number
of relevant crossing changes and an algorithm to find them all.

An algorithm to compute the constant $k$ is given in Section 11, building on Section 11 of [24], although it would
be technically challenging to implement.

Note that the finiteness statement in Theorem 1.2 does not obviously imply Theorem 1.1. Theorem 1.2
provides a bound on the number of ways of turning $L$ into a split link by a crossing change. But to
find this list of crossing changes is a highly non-trivial task. The algorithm that we give
is not very efficient, and we do not attempt
to provide an upper bound on its running time. Nevertheless, one can often apply
the techniques behind it quite practically. For example, we can obtain the following result.

\noindent {\bf Theorem 1.3.} {\sl Any crossing change that turns the Whitehead link
into a split link is equivalent to changing a crossing in some alternating diagram.}


An outline of the paper is as follows. In Section 2, we recall the operation of trivial tangle replacement.
This is a generalisation of a crossing change, and is in fact the central object of study within this paper.
It is well known that trivial tangle replacement can be studied by analysing the double cover of the 3-sphere
branched over the link, via the Montesinos trick. We recall the relevant theory in Section 2. In Section 3,
we compare trivial tangle replacement with crossing changes, focusing in particular on the notions
of equivalence in each case. In Section 4, we give a characterisation
of the hyperbolic links in the 3-sphere that are 2-string prime, in terms of their branched double
covers. We also consider double covers branched over sublinks of the link, which are the 3-manifolds
that play a central role in the proof of Theorem 1.1. In Section 5, we give an overview of the
general set-up of our algorithm. This involves constructing certain double branched covers $M$
and then searching for exceptional surgery curves in $M$. This second step uses earlier work of
the author [24]. The algorithm divides according to whether $M$ is Seifert fibred or hyperbolic.
In Section 6, the Seifert fibred case is analysed. In Section 7, the hypotheses of the main theorem in [24]
are verified. In Section 8, we analyse the mapping
class group of finite-volume hyperbolic 3-manifolds, mostly from an
algorithmic perspective.  In Section 9, we
show that all the problems in Theorem 1.1 are decidable.  In Section 10, we give an overview of the work in [24]. 
This leads to the finiteness result, Theorem 1.2, in Section 11. In Section 12, we analyse
the Whitehead link and we classify the crossing changes that
can be applied to the link to turn it into a split link. 

The author gratefully acknowledges the helpful suggestions of the referee, which have substantially improved this paper.

\vskip 18pt
\centerline {\caps 2. Tangle replacement}
\vskip 6pt

In this section, we recall the operation of tangle replacement, and the well-known Montesinos trick [35].

A {\sl tangle} is a 1-manifold $A$ properly embedded within a 3-ball $B$. When $A$ has no closed
components, and so is a collection of $k$ arcs for some positive integer $k$, it is termed a
{\sl $k$-string} tangle. The tangle is {\sl trivial}
if there is homeomorphism between $B$ and $D^2 \times I$ taking $A$ to $P \times I$, for some finite collection
of points $P$ in the interior of $D^2$. In the case of a trivial 2-string tangle, its {\sl core} is an arc
$\alpha \times \{ \ast \}$, where $\alpha$ is an embedded arc in the interior of $D^2$
joining the two points of $P$, and $\ast$ is a point in the interior of $I$. It is in fact the case
that a trivial 2-string tangle has a unique core up to isotopy of $B$ that leaves $A$ invariant.

Let $M$ be a compact orientable 3-manifold with (possibly empty) boundary. Let $L$ be a compact
1-manifold properly embedded in $M$. 
Let $\alpha$ be an arc embedded in the interior of $M$ such that $L \cap \alpha = \partial \alpha$.
Let $B$ be a regular neighbourhood of $\alpha$ in $M$, which intersects
$L$ in a trivial 2-string tangle in which $\alpha$ is a core arc. Suppose that we remove this tangle and insert
into $M$ another trivial tangle with the same endpoints. The result is
a new 1-manifold in $M$, which we say is obtained from $L$ by {\sl tangle replacement}
along $\alpha$. The possible trivial tangles that we may insert are
parametrised as follows. On $\partial B - L$, there is a unique isotopy class $\sigma$
of essential simple closed curves that bound a disc in the complement of the new tangle. We term
this the {\sl tangle slope}. The link that results from $L$ by this tangle replacement
is denoted $L_\sigma$.

We say that the {\sl distance} $\Delta(\sigma, \sigma')$ between two tangle slopes $\sigma$ and $\sigma'$ on $\partial B - L$ is
equal to half the minimal intersection number between two representative simple closed curves.
Any simple closed curve on $\partial B - L$ is separating, and hence any two curves have
even intersection number. Therefore, the distance between slopes is always an
integer.

Let $A$ be a trivial 2-string tangle in the 3-ball $B$. Then there is a unique double cover $V$ of $B$
branched over $A$. It is well known that $V$ is a solid torus. This is because $V$ is of the
form $A^2 \times I$, where $A^2$ is the annulus that is the double cover of the disc branched
over two points.

Consider the link $L_\sigma$ obtained from $L$ by tangle replacement.
Suppose that $M$ admits a double cover $\tilde M$ branched over $L$.
Then there is a corresponding double cover $\tilde M_\sigma$ branched over $L_\sigma$
which is defined as follows. The inverse images of $B$ and $M - {\rm int}(B)$ in $\tilde M$ give
double covers branched over $L \cap B$ and $L - {\rm int}(B)$ respectively.
Similarly, the inverse image of $\partial B$ is a double cover of the 2-sphere branched
over four points. Now there is a unique double cover of $S^2$ branched over four points.
This is a torus $T$.
Hence, any homeomorphism $\partial B \rightarrow \partial B$ that sends $\partial B \cap L$
to $\partial B \cap L$ lifts uniquely to a homeomorphism $T \rightarrow T$. One may view the tangle
replacement as simply attaching $B$ to $M - {\rm int}(B)$ via some homeomorphism that leaves
$\partial B \cap L$ invariant.
Lifting this homeomorphism to the branched double covers gives a gluing map, via
which we may construct the double cover $\tilde M_\sigma$ branched over $L_\sigma$.
Since the double cover of $B$ branched over $L \cap B$ is a solid torus, $\tilde M_\sigma$
and $\tilde M$ are related by Dehn surgery. The surgery curve in $\tilde M$ is the inverse image
of the arc $\alpha$. It is easy to check that the distance between 
between the two surgery slopes, one giving $\tilde M_\sigma$ and the other giving $\tilde M$,
is equal to the distance $\Delta(\sigma, \mu)$ between $\sigma$ and the meridian slope $\mu$ of $\alpha$.

The above use of branched double covers leads to a very useful method of parametrising slopes
of trivial tangles. We will consider trivial tangles $A$ within the 3-ball $B$, where $\partial B \cap A$ is
a given set of four points. The tangle is determined by the unique isotopy class of essential
curves in $\partial B - A$ that bound  a disc in the complement of $A$.
The 2-sphere $\partial B$ admits a unique double cover branched
over $\partial B \cap A$, which is a torus $T$. The elevation of the simple closed curve
in $\partial B - A$ that bounds a disc in $B - A$ is an essential simple
closed curve in $T$. One may parametrise the tangle $B \cap A$ by means of this slope.
It is possible to show (for example [5]) that this slope determines the trivial tangle up to an isotopy
of $B$ fixed on $\partial B$. Moreover, each slope
is realised by some tangle. Thus, trivial tangles are in one-one correspondence with 
slopes on $T$. One can pick a basis $\{ \lambda, \mu \}$ for the homology of $H_1(T)$,
and in the usual way, the slope with class $\pm (p \lambda + q \mu)$ in $H_1(T)$
is represented by $p/q \in {\Bbb Q} \cup \{ \infty \}$.
Given $p/q \in {\Bbb Q} \cup \{ \infty \}$, one may explicitly construct the associated tangle,
as follows. Let $[c_1, \dots, c_n]$ denote a continued fraction expansion for $p/q$, where
each $c_i \in {\Bbb Z}$. Then the associated tangle is shown in Figure 2, with the two possibilities
shown depending on whether $n$ is even or odd.

\vskip 12pt
\centerline{
\epsfxsize = 6in
\epsfbox{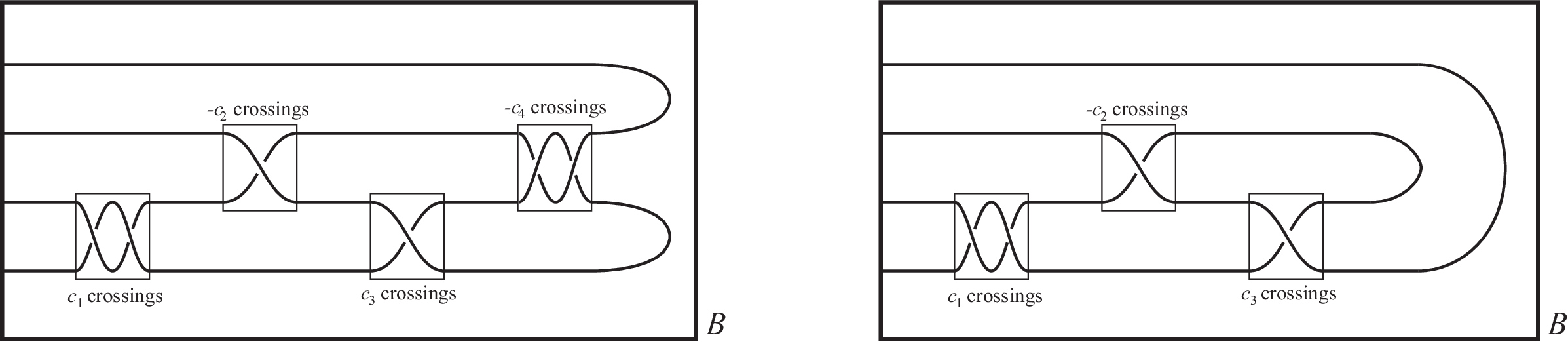}
}
\vskip 6pt
\centerline{Figure 2}

In the above figure, each box contains a line of crossings in a row, called a {\sl twist region}.
Conventionally, in a twist region with a `positive' number of crossings, they are twisted in 
a clockwise-fashion. However a box with a `negative' number $c$ of crossings
is in fact a string of $|c|$ crossings twisted in an anti-clockwise fashion.

We will be considering the operation of tangle replacement throughout this paper.
At various points, it will be important for us to consider when tangle replacement
can change a trivial tangle to another trivial tangle. More specifically,
suppose that $A$ is a trivial 2-string tangle in the 3-ball $B$, and suppose that $\alpha$ is an embedded
arc in the interior of $B$ such that $\alpha \cap A = \partial \alpha$. Suppose that tangle replacement along
$\alpha$ changes $A$ into another trivial 2-string tangle $A'$. Then what are the possible locations
for $\alpha$, what are the possible tangle replacements, and what is the relationship between the
slopes of $A$ and $A'$? Fortunately, all of these questions have been given a precise answer
by Baker and Buck [1], by use of branched double covers and surgical methods of Gabai [10].
The situation is simplest to state when the distance of the tangle replacement is at least two, as follows
(see Theorems 1.1 and 3.1 in [1]).

\noindent {\bf Theorem 2.1.} {\sl Let $A$ be a trivial 2-string tangle with slope $\infty$ in the 3-ball $B$.
Suppose that $\alpha$ is an embedded
arc in the interior of $B$ such that $\alpha \cap A = \partial \alpha$. Suppose that distance $d \geq 2$ tangle replacement along
$\alpha$ changes $A$ into another trivial 2-string tangle $A'$ with slope $p/q$. Then one of the
following holds:
\item{(i)} $\alpha$ is the core arc of the tangle $A$;
\item{(ii)} $p/q = (1 \pm dab)/\pm da^2$, for coprime integers $a$ and $b$. Moreover, if $a/b$ has continued fraction
expansion $[c_1, \dots, c_n]$, then there is an isotopy of $B$, fixed on $\partial B$, taking $B \cap A$
to the trivial tangle with continued fraction expansion $[0, c_1, \dots, c_n, 0, -c_n, \dots, -c_1]$, and taking $\alpha$
to the crossing arc of the central twist region labelled $0$. The tangle replacement simply
replaces the $0$ crossings with $\pm d$. (See Figure 3.)

}

\vskip 12pt
\centerline{
\epsfxsize = 5in
\epsfbox{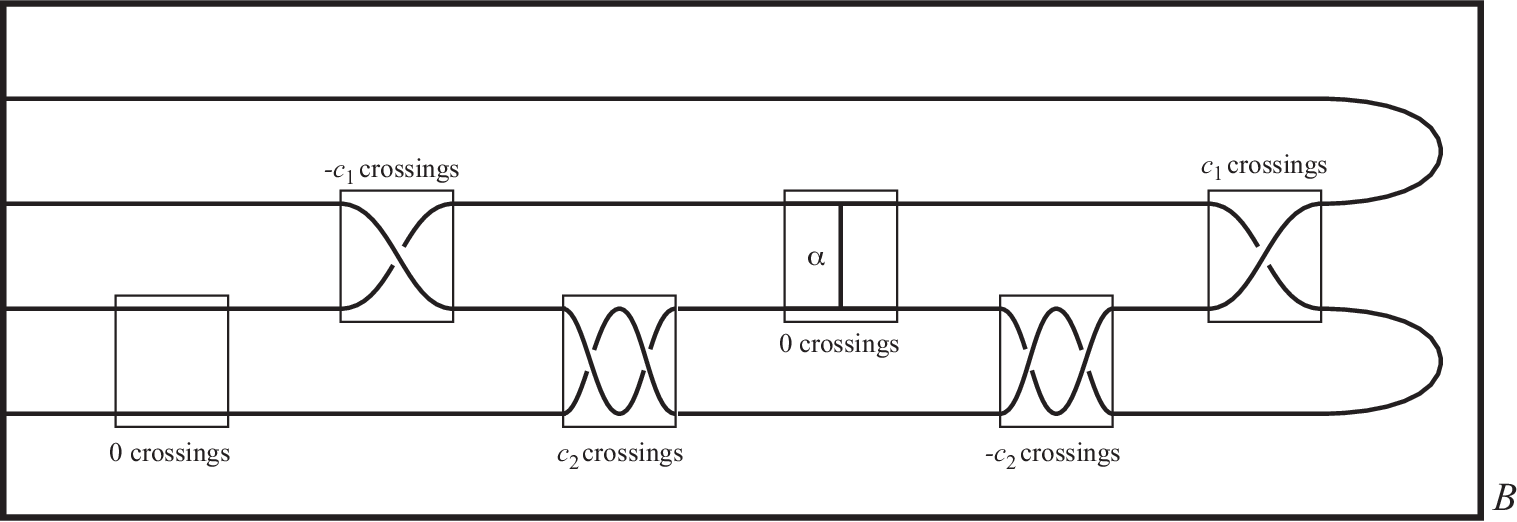}
}
\vskip 6pt
\centerline{Figure 3}

Using methods similar to those of Baker and Buck, we can obtain the following result.

\noindent {\bf Theorem 2.2.} {\sl Let $A$ be a 2-string tangle in the 3-ball $B$ such that $\partial B - \partial A$ is
compressible in the complement of $A$.
Let $\alpha$ be an embedded arc in the interior of $B$ such that $\alpha \cap A = \partial \alpha$. Suppose that
$\partial B - \partial A$ is incompressible in the complement of $A \cup \alpha$. Let $A'$ be obtained
from $A$ by trivial tangle replacement along $\alpha$ with distance at least two. Then the 
following hold.
\item{(i)} There is no 3-ball in $B$ with 
boundary disjoint from $A'$ and that encloses a closed component of $A'$.
\item{(ii)} If $\partial B - \partial A'$ is compressible in the complement of $A'$, then $A$ and $A'$ are trivial
tangles, and hence (i) or (ii) of Theorem 2.1 holds.

}

\noindent {\sl Proof.} Let $M$ and $M'$ be the double covers of $B$, branched over $A$ and $A'$
respectively. These differ by surgery along a curve $K$ in $M$ that is the inverse image of $\alpha$. The distance between the surgery
slope and the meridian slope is equal to the distance of the tangle replacement, which
is at least two by assumption. Note that $M$ has compressible boundary, since the inverse
image of a compression disc for $\partial B - \partial A$ contains a compression disc for $\partial M$.
Hence, $M$ is either a solid torus or reducible. In fact, $A$ is obtained from a trivial tangle by possibly tying
a little knot in one or both of its strings. Hence, $M$ is the connected sum of a solid torus with two rational
homology 3-spheres, one or both of which may be 3-spheres. In particular, $M$ contains no non-separating 2-sphere, and so the same is
true of $M - {\rm int}(N(K))$.

On the other hand, $\partial M$ is incompressible in 
the complement of $K$, for the following reason. If there were a compression disc for $\partial M$
in the complement of $K$, the equivariant disc theorem [33] would provide one or two disjoint compression
discs for $\partial M$ in the complement of $K$ that are invariant under the involution of $M$. These descend to a compression disc
$D$ for $\partial B - \partial A$ in the complement of $\alpha$ and that intersects $A$ in at most one point.
This disc $D$ cannot be disjoint from $A$ by hypothesis. Hence, it intersects $A$ in a single point. Its boundary
lies in the 2-sphere $\partial B$ and so bounds discs in $\partial B$. The union of either of these discs with $D$
forms a 2-sphere which, for parity reasons, must intersect $A$ an even number of times. Therefore,
$\partial D$ bounds a disc in $\partial B$ that intersects $\partial A$ once. The inverse image of $\partial D$
in $\partial M$ therefore bounds a disc in $\partial M$, which contradicts the fact that it is the boundary
of a compression disc.

To prove (i), suppose that $B$ contains a ball with boundary disjoint from $A'$ that encloses a closed component of $A'$.
The inverse image of this ball in $M'$ is a connected 3-manifold with two spherical boundary components.
The complement of this manifold in $M'$ is connected, and so we deduce that $M'$ contains
a non-separating 2-sphere. We now apply Scharlemann's theorem [45]. This implies that
a compact orientable irreducible 3-manifold with toroidal boundary cannot be Dehn filled along slopes
with distance at least two, and where one filling gives a manifold with compressible boundary and
the other filling gives a reducible 3-manifold. In our situation, $M - {\rm int}(N(K))$ need not be
irreducible, but it is a connected sum of irreducible 3-manifolds, one of which contains $\partial N(K)$.
Let $X$ be this summand. This summand must contain $\partial M$, as otherwise the meridional Dehn filling
of $M - {\rm int}(N(K))$ could not produce a 3-manifold with compressible boundary. So, $M - {\rm int}(N(K))$ is
the connected sum of $X$ with a rational homology 3-sphere (which may be a 3-sphere).
When $X$ is filled to give a summand of $M'$, this summand must contain a non-separating sphere
in $M'$. Hence, when $X$ is Dehn filled in two different ways with distance at least two, one filling gives a 3-manifold with
compressible boundary, and the other filling gives a reducible 3-manifold. This contradicts
Scharlemann's theorem.

We now prove (ii). We observed above that $M$ has compressible boundary.
Suppose that $\partial B - \partial A'$ is compressible in the complement of $A'$. Then
$M'$ also has compressible boundary. Theorem 2.4.4 of [8] therefore applies. With the
assumption that the distance of the surgery is at least two, it implies that $M - {\rm int}(N(K))$
is either a copy of $T^2 \times I$ or a `cable space'. In the former situation, every Dehn filling
of $K$ gives a solid torus. So, suppose that $M - {\rm int}(N(K))$ is a cable space. This is
a Seifert fibred space with annular base space and with one singular fibre. 
If one were to fill $\partial N(K)$ along a slope that has distance one from the regular fibre, the
result is a solid torus. In all fillings with distance at least two from the regular fibre, the resulting manifold has incompressible boundary,
because it is a Seifert fibre space with two singular fibres.
Since we are filling $\partial N(K)$ along slopes with distance at least two and we obtain manifolds
$M$ and $M'$ with compressible boundary, we deduce that the slopes giving $M$ and $M'$ have distance
$1$ from the regular fibre. Thus, $M$ and $M'$ are both solid tori. Since $M$ is
the branched double cover over $A$, it admits a (piecewise-linear) involution. Piecewise-linear involutions of the
solid torus have been classified (see Theorem 4.3 in [15]). Up to conjugacy by a piecewise-linear homeomorphism,
there is just one orientation-preserving piecewise-linear involution with fixed-point set homeomorphic to two intervals.
Therefore, $A$ is a trivial tangle. Similarly, $A'$ is a trivial tangle. We are therefore in the setting of Theorem 2.1,
and hence (i) or (ii) of Theorem 2.1 holds. $\square$

\vskip 18pt
\centerline{\caps 3. Crossing arcs versus crossing circles}
\vskip 6pt

A crossing change to a link $L$ can be viewed in two ways: as a special type of tangle replacement
and as a special type of Dehn surgery. In this section, we will explore how these two alternative viewpoints
are related.

Tangle replacement was discussed in Section 2. One starts with an embedded arc $\alpha$ such that
$\alpha \cap L = \partial \alpha$. A regular neighbourhood of $\alpha$ in $S^3$ intersects $L$ in a trivial
tangle. The crossing change is implemented by removing this tangle and inserting another trivial
tangle, with the property that the new and the old tangle slopes have distance exactly 2. There are
infinitely many possible tangle replacements of this form, as shown in Figure 4. Each corresponds to changing
a crossing in some diagram of $L$.

\vskip 12pt
\centerline{
\epsfxsize = 3.5in
\epsfbox{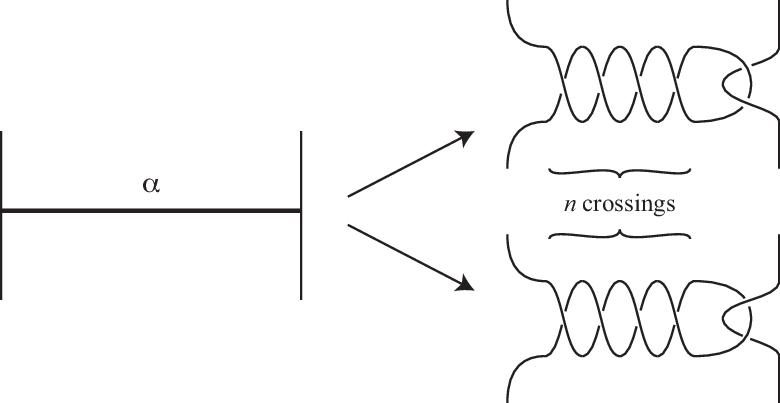}
}
\vskip 6pt
\centerline{Figure 4}
\vfill\eject

This ambiguity is a slightly unfortunate one. It could, of course, be rectified by requiring $\alpha$ to be framed in some way.
More precisely, one could specify not just the arc $\alpha$ but also an explicit identification between its regular neighbourhood $B$
and $D^2 \times I$, so that $B \cap L$ is sent to vertical arcs in the product structure. One would then be able
to specify the precise tangle replacement by giving the slope of the new tangle as an explicit fraction.

This is somewhat cumbersome and so it is more usual to specify crossing changes via surgery along crossing circles,
as described in Section 1. In this section, we investigate the following questions. If a crossing change is
specified by tangle replacement along an arc, then how many crossing circles does this give rise to?
If a crossing change is specified by surgery along a crossing circle, how many associated
crossing arcs are there? 

The second of the above questions has a possibly surprising answer. A crossing circle can give rise to
an arbitrarily large number of distinct crossing arcs. The point is that, to obtain a crossing arc from a crossing circle,
one must choose a crossing disc $D$. The associated crossing arc is then the embedded arc in $D$ joining the two
points of $D \cap L$. But a crossing circle may bound many, quite different crossing discs. An example
is given in Figure 5, where a single crossing circle gives rise to 4 different crossing arcs. This example generalises
in an obvious way to arbitrarily many different crossing arcs.

\vskip 6pt
\centerline{
\epsfxsize = 2.5in
\epsfbox{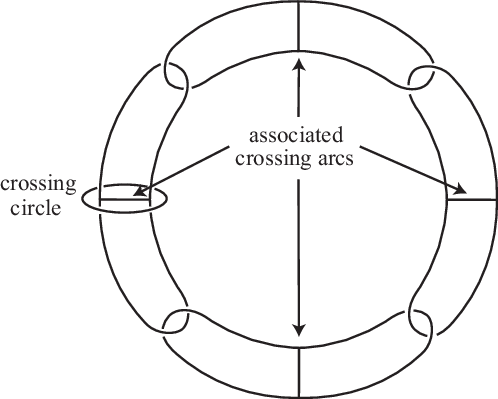}
}
\vskip 6pt
\centerline{Figure 5}

Let us now pass to the first of the above questions. When a crossing change is specified by a tangle replacement along $\alpha$, then there
are actually {\sl two} associated crossing circles, as shown in Figure 6.

\vskip 12pt
\centerline{
\epsfxsize = 3.5in
\epsfbox{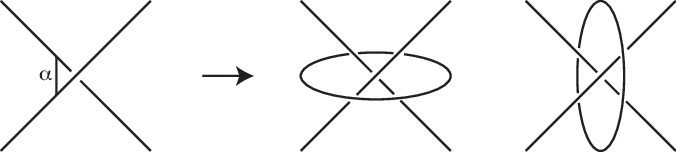}
}
\vskip 6pt
\centerline{Figure 6}

Implicit in the above statement is that no other crossing circles arise from this
tangle replacement. We now make this more precise.

\noindent {\bf Lemma 3.1.} {\sl When a crossing change to a link $L$ is achieved by tangle replacement
along an arc $\alpha$, this gives rise to precisely two crossing circles, so that surgery
along either of these crossing circles implements this crossing change.}

\noindent {\sl Proof.} Note first that the two crossing circles shown in Figure 6 are not equivalent. In other words,
they are not related by an ambient isotopy that preserves the link $L$. This is because one can pick any
orientation on $L$ and then one of these crossing circles has zero linking number with $L$ and other has linking 
number $\pm 2$.

We now need to show that we do not obtain any further crossing circles. To create a crossing circle from the crossing
arc $\alpha$, we must thicken $\alpha$ to a disc. There are infinitely many ways of doing this, that are parametrised
by an integer $n \in {\Bbb Z}$. We denote the boundary of this disc by $C_n$, as shown in Figure 7.

\vskip 12pt
\centerline{
\epsfxsize = 5in
\epsfbox{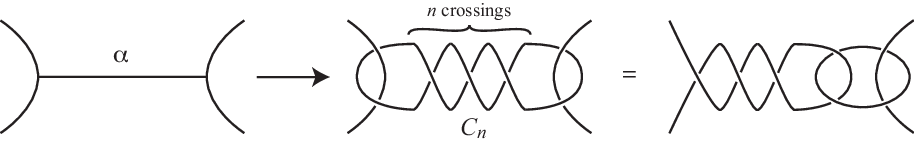}
}
\vskip 6pt
\centerline{Figure 7}

When $\pm 1$ surgery is performed along $C_n$, we obtain the tangle with slope $n \pm {1 \over 2}$,
according to Figure 2. Thus, we see that, to obtain a specific slope, in other words to obtain a specific
tangle replacement, there are exactly two choices of $n$. Specifically, we could perform
$-1$ surgery along $C_{n+1}$ or $+1$ surgery along $C_n$. $\square$

\vskip 18pt
\centerline{\caps 4. Simplicity of the branched double covers}
\vskip 6pt

Theorem 1.1 only applies to links that are hyperbolic and 2-string prime. It is reasonable to ask whether
it is possible to easily determine whether this condition holds. In the following result, we give
an alternative characterisation in terms of the geometry of the branched double cover.
This is easily checked in practice (using Snappea for example [52]) and can be determined
algorithmically [19, 22, 28, 46]. We also examine branched double covers over sublinks of the link $L$,
and derive a result that will be useful in the proof of Theorem 1.1.

\noindent {\bf Proposition 4.1.} {\sl Let $L$ be a hyperbolic link in the 3-sphere. 
\item{(i)} Then $L$ is 2-string prime if and only if
its branched double cover is hyperbolic or a small Seifert fibre space.
\item{(ii)} Suppose that $L$ is 2-string prime, and let $L'$ be
a sublink of $L$. Then the double cover of $S^3 - {\rm int}(N(L - L'))$
branched over $L'$ is hyperbolic or a small Seifert fibre space.

}

Recall that a Seifert fibre space is {\sl small} if it contains no essential embedded torus. In particular, the 3-sphere is a small
Seifert fibre space, as is any lens space.

\noindent {\sl Proof.} Note that the forwards implication in (i) is a special case of (ii) with $L' = L$.
So, we initially focus on (ii). Let $M$ be the double cover of $S^3 - {\rm int}(N(L - L'))$
branched over $L'$. To verify that $M$ is hyperbolic or a small Seifert fibre space, we appeal to the solution to 
the Geometrisation Conjecture [40, 41, 42]. So, if $M$ is not hyperbolic or a small Seifert fibre space,
then it is toroidal or reducible. Suppose first that $M$ is reducible. Then the equivariant sphere theorem (see Theorem 3 in [32] and its proof) 
implies that there are one or two embedded disjoint essential spheres that are invariant under the covering involution. 
Their union descends to a 2-sphere in $S^3$  either that is disjoint
from $L$ or that intersects $L$ in two points. If the sphere is disjoint from $L$, then it has components of $L$ on both sides, since 
its inverse image in $M$ is essential. If the sphere intersects $L$ in two points, then this forms an essential annulus properly
embedded in the exterior of $L$. In both cases, $L$ fails to be hyperbolic. Suppose now that $M$ is toroidal. Then the 
equivariant torus theorem (Corollary 4.6 in [15]) gives one or two embedded disjoint essential tori that are invariant under the covering involution. These descend to
an essential torus in the exterior of $L$ or to an essential 4-times punctured sphere with meridional boundary.
In the former case, this implies that $L$ is not hyperbolic. In the latter case, $L$ is not 2-string prime.

We now prove the backwards implication in (i). Let $M$ be the double cover of $S^3$ branched over $L$.
Suppose that $M$ is hyperbolic or a small Seifert fibre space. Let $S$ be a 2-sphere in $S^3$
that intersects $L$ in four points, such that $S - L$ has no compression disc in the complement of $L$. 
The inverse image of this 2-sphere in $M$ is a torus $T$, which must be compressible by our hypothesis about $M$. 
By the equivariant disc theorem (Theorem 7 in [33]), there
are one or two disjoint compression discs for $T$ that are invariant under the involution
of $M$. These descend to a compression disc $D$ for $S$ that intersects $L$ in at most one point.
This cannot be disjoint from $L$ by our assumption about $S$. On the other hand, if $D$ intersects $L$ in a single point, 
then $\partial D$ separates $S$ into two discs, one of which contains a single point of $L \cap S$.
We deduce in that case that the inverse image of $D$ was not a compression disc
for $T$, which is a contradiction. $\square$

\vskip 18pt
\centerline{\caps 5. The general set-up}
\vskip 6pt

Let $L$ be our given link in $S^3$. Suppose that a crossing change to $L$ transforms it into
a split link $L_\circ$. Associated with this crossing change is a crossing circle $C$ in the
complement of $L$ that bounds a crossing disc $D$. Running between the two points of $L \cap D$
is the crossing arc $\alpha$.

Let $L'$ be the union of the components of $L$ containing $L \cap D$. Thus, $L'$
has one or two components. We let $M$ be the manifold obtained from 
$S^3 - {\rm int}(N(L - L'))$ by taking the double cover branched over $L'$.
More precisely, we consider the double cover of $S^3 - {\rm int}(N(L))$ determined
by the homomorphism $\pi_1(S^3 - {\rm int}(N(L)) \rightarrow {\Bbb Z}/2$ that
measures the mod $2$ linking number of a loop with $L'$. Then $M$ is obtained
from this cover by Dehn filling each component of the inverse image of
$\partial N(L')$ using slopes that are elevations of meridians. 

Suppose that $L_\circ$ is obtained from $L$ by performing surgery along $C$
via the slope $\pm 1$. Let $L'_\circ \subset L_\circ$ be the image of $L'$ after this surgery.
Let $M_\circ$ be the double cover of $S^3 - {\rm int}(N(L_\circ - L'_\circ))$ branched over $L'_\circ$.
Then $M_\circ$ is obtained from $M$
by Dehn surgery along a curve $K$. Moreover, if $\mu$ is the meridional slope
on $\partial N(K)$ and $\sigma$ is the surgery slope, then $\Delta(\sigma, \mu) = 2$.

There are two main reasons why we use this set-up involving branched double covers,
rather than simply considering surgery along the crossing circle. Firstly, the distance between the
surgery slope $\sigma$ and the meridian slope $\mu$ is more than $1$. Secondly, we will
see that $H_2(M - {\rm int}(N(K)), \partial M) \not= 0 $. These two seemingly technical
points are important because Dehn surgery theory works most smoothly when they hold.
In particular, they are hypotheses in the following theorem of the author [24].

Let $M$ be a compact orientable 3-manifold with $\partial M$ a (possibly empty) union of tori.
Let $K$ be a knot in $M$, and let $\sigma$ be a slope on $\partial N(K)$ other than the
meridional slope $\mu$. Let $M_K(\sigma)$ be the manifold that is obtained by
Dehn surgery along $K$ via the slope $\sigma$. Then $\sigma$ is an {\sl exceptional slope}
and $K$ is an {\sl exceptional surgery curve} if any of the following holds:
\item{(i)} $M_K(\sigma)$ is reducible,
\item{(ii)} $M_K(\sigma)$ is a solid torus, or
\item{(iii)} the core of the surgery solid torus has finite order in
$\pi_1(M_K(\sigma))$.

\noindent Also, $\sigma$ and $K$ are {\sl norm-exceptional} if there is some
$z \in H_2(M - {\rm int}(N(K)), \partial M)$ that maps to an element 
$z_\sigma \in H_2(M_K(\sigma), \partial M_K(\sigma))$, such that the
Thurston norm of $z_\sigma$ is less than the Thurston norm of $z$.

\noindent {\bf Theorem 5.1.} {\sl There is an algorithm that takes, as its input,
a triangulation of a compact connected orientable 3-manifold $M$,
with $\partial M$ a (possibly empty) union of tori. The output to the algorithm
is a list of all knots $K$ within $M$ and all slopes $\sigma$ on $\partial N(K)$ with all the following properties:
\item{(i)} $M - {\rm int}(N(K))$ is irreducible and atoroidal, and $H_2(M - {\rm int}(N(K)), \partial M) \not= 0$;
\item{(ii)} $\sigma$ is an exceptional or norm-exceptional slope on $\partial N(K)$, such that
$\Delta(\sigma, \mu) > 1$, where $\mu$ is the meridian slope on $\partial N(K)$.

\noindent In particular, there are only finitely many such knots $K$ and slopes $\sigma$.}

The way that the algorithm lists the possibilities for $K$ is described in Section 10.
It is straightforward to then realise each possibility for $K$ as a subcomplex of a
suitable iterated barycentric subdivision of the given triangulation of $M$.
(See Theorem 10.2 and the discussion after it.)

Note that in our setting, $\sigma$ is an exceptional surgery slope on $\partial N(K)$.
This is because Dehn filling $M - {\rm int}(N(K))$ along $\sigma$ gives the manifold
$M_\circ$. This is the branched double cover of $S^3 - {\rm int}(N(L_\circ - L_\circ'))$
branched over $L'_\circ$. The splitting sphere in the complement of $L_\circ$ lifts to
reducing spheres in $M_\circ$.

Thus, roughly speaking, the algorithm required by Theorem 1.1 proceeds by constructing the finitely many possibilities for $M$, then using
Theorem 5.1 to find all the exceptional surgery curves $K$ in $M$ satisfying the hypotheses of Theorem 5.1, and then determining
whether any of these descend to a crossing arc $\alpha$ for $L$. It thereby
builds a finite list of possibilities for $\alpha$ and for each such possibility, it
provides the associated tangle replacement slope. In each case, we perform
this tangle replacement and determine whether the resulting link is split, totally
split or the unlink, as appropriate.

\vskip 18pt
\centerline {\caps 6. The Seifert fibred case}
\vskip 6pt

As in the previous section, $M$ is the double cover of $S^3 - {\rm int}(N(L - L'))$ branched over $L'$.
We saw in Proposition 4.1 that the manifold $M$ is hyperbolic or a small Seifert fibre space.
In this section, we deal with the case where $M$ is Seifert fibred. In this setting, the list
of potential exceptional surgery curves is very simple, as given by the following result.

\noindent {\bf Theorem 6.1.} {\sl Let $M$ be a Seifert fibred 3-manifold with non-empty boundary and let $K$ be
a knot in $M$ such that $M - {\rm int}(N(K))$ is irreducible and atoroidal and $H_2(M - {\rm int}(N(K)),\partial M) \not= 0$.
Suppose that $\sigma$ is an exceptional slope on $\partial N(K)$ such that $\Delta(\sigma, \mu) > 1$,
where $\mu$ is the meridional slope. Then $K$ is isotopic to a singular fibre of $M$ and $\sigma$ is
the slope of the regular fibres, when $N(K)$ is a fibred regular neighbourhood of $K$.
In particular, $M - {\rm int}(N(K))$ is Seifert fibred.}

We will defer the proof of this until Section 10.

The above result says that the exceptional surgery curve $K$ is isotopic to a singular
fibre of $M$. However, the manifold $M$ that we are considering comes with an involution
that preserves $K$, and there is no {\sl a priori} reason why this isotopy should be
equivariant with respect to the involution. We deal with this as follows.

\noindent {\bf Addendum 6.2.} {\sl Let $M$, $K$, $\sigma$ and $\mu$ be as in Theorem 6.1.
Let $\tau$ be a piecewise-linear involution of $M$ that leaves $K$ invariant. Then there is a Seifert fibration of
$M$ that is invariant under $\tau$ and that has $K$ as a singular fibre. Again, $\sigma$ is the 
slope of the regular fibres when $N(K)$ is a fibred regular neighbourhood of $K$.}

One might wonder why we require that $\tau$ is piecewise-linear. Indeed, Theorem 6.1 does not refer
to any triangulation of $M$. But by [34], $M$ has a unique piecewise-linear structure, and so the requirement that
$\tau$ is piecewise-linear is a well-defined condition. The reason we make this assumption is that an involution
$\tau$ of a 3-manifold $M$ need not be conjugate to a piecewise-linear involution and hence $M / \tau$ need not
even be triangulable [2]. We impose the hypothesis that $\tau$ is piecewise-linear to avoid this
situation.

Before we prove Addendum 6.2, we quote the following theorem.

\noindent {\bf Theorem 6.3.} {\sl Let $M$ be a Seifert fibre space that admits a piecewise-linear involution
$\tau$. Then it admits a Seifert fibration that is invariant under $\tau$.}

This was proved by Tollefson [51], but under the assumption that if the base space of
$M$ is a 2-sphere, then it has at least four singular fibres. The case excluded by Tollefson requires
the Orbifold Theorem [3, 7]. The quotient $M / \tau$ is an orbifold. Since $M$ is Seifert fibred,
and its base space is a 2-sphere with at most $3$ singular points, then $M/  \tau$ is orbifold-irreducible
and orbifold-atoroidal. Hence by the Orbifold Theorem, it is either hyperbolic or Seifert fibred.
But if $M / \tau$ were hyperbolic, then so would $M$ be, which is impossible. Thus,
$M / \tau$ is Seifert fibred. Its Seifert fibration lifts to a Seifert fibration of $M$ that is
invariant under $\tau$.

\noindent {\sl Proof of Addendum 6.2.} We assume Theorem 6.1, which will be proved in Section 10.
By Theorem 6.1, $M - {\rm int}(N(K))$ is Seifert fibred. It admits an involution, which is the restriction
of $\tau$. Hence, by Theorem 6.3, it admits a Seifert fibration that is invariant under the involution.
The slope of the regular fibres must be the exceptional slope $\sigma$, because filling
along any other slope gives a Seifert fibre space. Moreover, this Seifert fibre space cannot be a
solid torus, by our hypothesis that $H_2(M - {\rm int}(N(K)),\partial M) \not= 0$.
Thus, the Seifert fibration on $M - {\rm int}(N(K))$ extends to a Seifert
fibration on $M$ that is invariant under the involution $\tau$. Since the meridional slope $\mu$
is assumed to have distance at least two from the slope of the regular fibre $\sigma$, we
deduce that the surgery curve $K$ is a singular fibre. $\square$

\noindent {\bf Remark 6.4.} Note that Addendum 6.2 and Theorem 6.3 assert the existence of
{\sl some} Seifert fibration that is invariant under $\tau$. However, most Seifert fibre spaces have
a unique Seifert fibration up to isotopy. Indeed, this is true of any Seifert fibre space with non-empty
boundary other than $T^2 \times I$, the twisted $I$-bundle over the Klein bottle and the solid torus
(see Theorem VI.18 of [16] for example). This is important, because once one has identified
a Seifert fibration on the manifold $M$ in Theorem 6.1, then we know that $K$ is isotopic
to a singular fibre in this Seifert fibration. In the setting of Addendum 6.2, $K$ is also
assumed to be invariant under the piecewise-linear involution $\tau$, and hence it descends
to a 1-manifold in the orbifold $M/ \tau$. We would like to know that $K$ is similarly
well-defined up to isotopy in $M/\tau$. Fortunately, this is provided to us by the following
result of Bonahon and Siebenmann (see Theorem 2 in [4]) which establishes the uniqueness
up to isotopy of Seifert fibrations on many orbifolds.

\noindent {\bf Theorem 6.5.} {\sl Let $M / \tau$ be a Seifert fibred orbifold with non-empty boundary. Suppose that
the base orbifold of the Seifert fibration is not finitely covered by a disc or an annulus. Then the Seifert
fibration on $M / \tau$ is unique up to isotopy that preserves the singular locus of the orbifold throughout.}

We will need the following constructive version of Theorem 6.3.

\noindent {\bf Theorem 6.6.} {\sl There is an algorithm that takes, as its input, a triangulation
$T$ of a Seifert fibre space $M$ with non-empty boundary and an involution $\tau$ of $M$ that
preserves $T$. The output of the algorithm is a union of disjoint simple closed curves that 
are the singular fibres in a Seifert fibration of $M$ that is invariant under $\tau$.
The algorithm also produces a regular neighbourhood of these singular fibres, together
with a slope on each of these solid tori that represents a regular fibre.}

\noindent {\sl Proof.} In Algorithm 8.1 of [19], Jaco and Tollefson provided an algorithm to determine
whether a Haken manifold with incompressible boundary is Seifert fibred. It also produced the information required by the theorem:
the singular fibres, together with the slopes of regular fibres on the boundary of their solid toral neighbourhoods.
However, we need to perform a version of this procedure equivariantly.

If necessary, we first subdivide the triangulation $T$ to a triangulation $T'$ so that the fixed-point set of $\tau$ is 
simplicial in $T'$.

If $M$ is a solid torus, then this may be algorithmically determined, for example using Theorem 6.2 in [19]. In this case, $M$
admits a Seifert fibration that is invariant under $\tau$. It has at most one singular fibre, that is a core of $M$. Hence,
our algorithm must simply find the slope of the regular fibres on $\partial M$. But this slope may be taken to be
any non-meridional slope that is invariant under the involution, and this may easily be determined algorithmically.

If $M$ is homeomorphic to $T^2 \times I$, this may also be algorithmically determined, using Algorithm 8.1 in [19].
In this case, our algorithm ends by declaring that $M$ has no singular fibres.

So, we may assume that $M$ is neither a solid torus nor $T^2 \times I$. It therefore contains an essential properly embedded annulus. 
According a theorem of Kobayashi (Theorem 1 in [21]), it contains such an annulus that is either invariant under $\tau$
or disjoint from its image under $\tau$. We need to show that in addition, this annulus $A$
can be realised as a normal surface, with control over its number of intersections with the 1-skeleton of $T'$.

By the PL-minimal surface theory of Jaco and Rubinstein, there is a PL-least area surface in the isotopy class of $A$,
which we will also call $A$. By definition, this is normal in $T'$. By Theorem 7 of [18], $A$ is either disjoint
from its image under $\tau$ or it equals its image. Let $\tilde A$ be $A \cup \tau A$, which is equal either to 
$A$ or to the disjoint union of $A$ and $\tau A$.

Normal surface theory [14, 17, 30] gives that there is a finite constructible collection of normal surfaces $F_1, \dots, F_n$
in $T'$ that are fundamental. The normal surface $\tilde A$ is a normal sum $k_1 F_1 + \dots + k_n F_n$. We will show that 
$\tilde A$ can be chosen so that each $k_i$ is at most $8$. By Theorem 4.1.36 of [30] (see also Theorem 2.2 of [17]), 
any normal summand of $\tilde A$ must be incompressible and boundary-incompressible, and no summand
can be a sphere or disc. Since Euler characteristic is additive under normal summation, 
any fundamental surface $F_i$ that is a summand of $\tilde A$ must be an annulus or M\"obius band.
If $F_i$ is a M\"obius band, then $2F_i$ is an annulus. Hence if some $k_i > 1$, then $\tilde A$
has an annulus $A'$ as a summand which is either fundamental or twice a fundamental surface.
The surface $\tau A'$ is also a normal surface. It is also a summand for $\tilde A$.
Hence, $A'$ and $\tau A'$ have compatible normal co-ordinates in the sense that no tetrahedron of
$T'$ contains a quadrilateral of $A'$ and a quadrilateral of $\tau A'$ that are not normally isotopic.
One may therefore form the normal sum $A' + \tau A'$ and obtain an embedded normal surface. 
Now $A' + \tau A'$ is a summand of $\tilde A + \tau \tilde A = 2\tilde A$. Therefore
by Theorem 4.1.36 in [30], no component of the surface $A' + \tau A'$ is a sphere or disc. Hence, it 
is a union of annuli and M\"obius bands. By Theorem 4.1.36 in [30], these are incompressible and
boundary-incompressible. Note that $A' + \tau A'$ is invariant under $\tau$ up to normal
isotopy. In fact, $A' + \tau A'$ has a normal representative that is actually invariant under $\tau$ by Theorem 2 in [18].
Pick a component of $A' + \tau A'$. It is either invariant under $\tau$ or disjoint from its image.
If this is an annulus, then it is the required surface, because it is a sum of at most $4$ fundamental surfaces. 
On the other hand, if this component of $A' + \tau A'$ is a M\"obius band, then its normal sum with itself
is the required annulus.

Thus, by searching through normal surfaces of the form $k_1 F_1 + \dots + k_n F_n$, where
each $k_i \leq 8$, we eventually find a normal essential annulus $A$ that is either invariant
under $\tau$ or disjoint from its image under $\tau$. We now subdivide the triangulation $T'$ equivariantly so that $A \cup \tau A$ is a subcomplex of it.
We can cut along $A \cup \tau A$, to form a triangulation of a new Seifert fibred manifold with toral boundary components.
Repeating in this way, we eventually decompose our Seifert space into a union of fibred solid tori.
In their boundary are a collection of annuli, which are copies of the last annulus or annuli that we
decomposed along to form the relevant solid torus. The slopes of these annuli are the slopes of the regular
fibres. We can thereby determine, for each of these solid tori, whether they have a singular fibre as
a core curve. Our algorithm ends by outputting these solid tori, their core curves and the slopes
of the regular fibres on their boundary. $\square$

\vskip 18pt
\centerline{\caps 7. Verifying the hypotheses of Theorem 5.1.}
\vskip 6pt

We now spend some time verifying that the hypotheses of Theorem 5.1
do hold in our setting.

We need to verify that
$H_2(M - {\rm int}(N(K)),\partial M) \not= 0$. Now it is not hard to check that
$H_2(M - {\rm int}(N(K)), \partial M)$ is a subgroup of $H_2(M, \partial M)$, and
its rank is either equal to that of $H_2(M, \partial M)$ or one less. By Poincar\'e
duality, the rank of $H_2(M, \partial M)$ is equal to the first Betti number of $M$,
and this is at least the number of toral boundary components. Each component
of $L - L'$ gives rise to one or two components of $\partial M$. So, when $|L| \geq 4$,
then $|\partial M| \geq 2$ and so we deduce that $H_2(M - {\rm int}(N(K)),\partial M) \not= 0$.

When $|L| = 2$, we are assuming that the two components of $L$ have zero
linking number. So, the crossing disc $D$ that $C$ bounds must intersect a single component
of $L$. For otherwise, the crossing change would modify the linking number by $\pm 1$,
and the result could not be a split link. So, $|L'| = 1$. Moreover,
the linking number between $L'$ and $L - L'$ is zero, and so
$M$ has two boundary components. Once again we deduce that
$H_2(M - {\rm int}(N(K)),\partial M) \not= 0$. 

The final case is where $|L| = 3$. This is a little more delicate. The disc $D$ can
intersect at most two components of $L$, and so there is a component $L_1$
that is disjoint from $D$. Let $\alpha$ be the crossing arc in $D$ running between the
two points of $D \cap L$. Let $S$ be a Seifert surface for $L_1$.
Since $\alpha$ is an arc, we may slide any points of $\alpha \cap S$
along $\alpha$ and off it. This may introduce new points of intersection between
$S$ and $L - L_1$ but $S$ remains a Seifert surface for $L_1$, and
so we may assume that $S$ is disjoint from
$\alpha$. The inverse image of $\alpha$ in $M$ is the surgery curve $K$.
The inverse image of $S$ in $M$ is a properly embedded orientable non-separating surface.
This represents a non-trivial element of $H_2(M - {\rm int}(N(K)),\partial M)$,
which verifies that this group is non-zero.

There are two remaining hypotheses in Theorem 5.1:
that $M - {\rm int}(N(K))$ is irreducible and atoroidal. If $M - {\rm int}(N(K))$ is
reducible, then the equivariant sphere theorem states that there are one or two disjoint reducing 2-spheres
in $M - {\rm int}(N(K))$ that are invariant under the involution. They descend
to a 2-sphere in $S^3 - {\rm int}(N(L - L'))$ that intersects $L'$ in either two or zero points.
Hence, by the hypothesis that $L$ is hyperbolic, this sphere $S$ bounds a ball $B$
so that $B \cap L$ is either empty or a trivial 1-string tangle. The
sphere or spheres lie in $M - {\rm int}(N(K))$, and so their image $S$ is disjoint from $N(\alpha)$.
Therefore, $\alpha$ must lie in $B$, because otherwise each component of the inverse image of $B$
is a ball in $M - {\rm int}(N(K))$. The
tangle replacement occurs within $B$, and so $B \cap L$ is replaced by a possibly non-trivial 1-string tangle.
However, this cannot make the link $L_\circ$ split, which is a contradiction.

Note that here, we used the fact that we are performing a crossing change to $L$. Thus, this
argument does not immediately extend to other tangle replacements. However, there
is another argument that works in this more general setting. Suppose that tangle
replacement is performed along $\alpha$ and that this turns the trivial 1-string tangle $B \cap L$
into a tangle that is split. Then on passing to the branched double cover, we deduce that
surgery along on a knot in the 3-ball creates a manifold containing a non-separating sphere.
Hence, the distance between the surgery slope and the meridian slope is $1$.
So, if we assume that the distance is more than $1$, then we reach a contradiction.
In fact, by using Gabai's proof of the Property R conjecture [9] and the solution to
the Smith conjecture [36], we would be able to classify the possible tangle replacements
even in the distance 1 case.

Suppose now that $M - {\rm int}(N(K))$ is toroidal. The equivariant torus theorem
gives that there are one or two disjoint essential embedded tori in $M - {\rm int}(N(K))$ that are invariant under the involution. They descend
either to a sphere in $S^3 - {\rm int}(N(L - L'))$ that intersects $L'$ in four points, or to a
torus disjoint from $L$. We consider these two cases separately.

Suppose that the torus or tori in $M - {\rm int}(N(K))$ project to an embedded torus $T$ in
$S^3 - {\rm int}(N(L))$. We are assuming that $L$ is a hyperbolic link, and so $T$ must bound a 
solid torus in the complement of $L$, or must lie within 
a 3-ball in the complement of $L$, or must be parallel
to a component of $\partial N(L)$. 

Suppose that $T$ bounds a solid torus in the
complement of $L$. Since its inverse image in $M$ is disjoint from $K$, $T$ is disjoint from $\alpha$.
As $\alpha$ starts and ends on $L$, it is therefore disjoint from the solid torus. Hence,
the torus was not essential in $M - {\rm int}(N(K))$, which is a contradiction. 

Suppose that
$T$ lies within a 3-ball in the complement of $L$ but does not bound a solid torus in the complement of $L$.
Then $T$ separates $S^3$ into two components, one of which is disjoint from $L$ and is homeomorphic to the
exterior of a non-trivial knot. The other component must be a solid torus containing $L$. As $\alpha$ has its endpoints on
$L$ and is disjoint from $T$, we deduce that $\alpha$ lies in this solid torus.
A meridian disc for this solid torus is disjoint from $L$, because the torus lies in a 3-ball disjoint from $L$.
Hence, every curve on $T$ has zero linking number with every component of $L$. Therefore,
the inverse image of $T$ in the branched double cover $M$ is two tori. Together these bound the inverse
image $V$ of the solid torus. The surgery curve $K$ lies in $V$, since $K$ is the inverse image of $\alpha$.
Note that $V - {\rm int}(N(K))$ is irreducible as otherwise the equivariant sphere theorem implies
that there is an essential sphere in $S^3 - {\rm int}(N(L \cup \alpha))$, which would imply that $L$ is split,
contrary to assumption. Since $T$ is compressible in $S^3 - {\rm int}(N(L))$, $\partial V$ is compressible in $V$. Since
$V$ has more than one toral boundary component, it is therefore reducible.
When surgery on $M$ is performed along $K$, a reducible manifold $M_\circ$ is created.
Let $V_\circ$ be the submanifold of $M_\circ$ that comes from $V$; so $V_\circ$ is obtained from $V$ by
surgery along $K$. Since $M_\circ$ is reducible, we deduce that $V_\circ$ is reducible
or has compressible boundary. Again because $V_\circ$ has more than one toral boundary
component, it must be reducible. A theorem of Gordon and Luecke [12] then states that
when an irreducible 3-manifold (in this case, $V - {\rm int}(N(K))$) with a toral boundary component
is Dehn filled in two different ways to obtain reducible 3-manifolds (in this case, $V$ and $V_\circ$),
the distance between the surgery slopes is one. This contradicts our assumption that
the distance of the tangle replacement is at least two.

Finally consider the case where $T$ is parallel to a component of $\partial N(L)$.
Then $T$ bounds a solid torus $W$ in $S^3$, that intersects $L$ in a single core curve.
Since $T$ is the image of a torus disjoint from $K$, $T$ is disjoint from $\alpha$.
If $\alpha$ is disjoint from the solid torus bounded by $T$, then the inverse image
of $T$ in $M - {\rm int}(N(K))$ is boundary parallel, which is contrary to hypothesis.
So, $\alpha$ lies in the solid torus bounded by $T$. Now, $T$ does not admit a
compression disc in $S^3 - {\rm int}(W)$ that is disjoint from $L$. This is because $T$
would then be compressible in the complement of $L \cup \alpha$, and hence the inverse
image of $T$ in $M - {\rm int}(N(K))$ would be compressible. After the tangle replacement,
$T$ continues to bound a solid torus in $S^3$, but its intersection with the new link need not be
a core curve. However, it still has winding number one, and so $T$ remains incompressible
in the complement of the new link $L_\circ$. Now this link complement contains an essential sphere,
because $L_\circ$ is split. After modifying this 2-sphere appropriately, we can make 
it disjoint from the incompressible torus. It then is disjoint from the solid torus,
and hence is disjoint from the inserted tangle. It therefore corresponds to a 2-sphere
in the complement of $L$ that separates components of $L$. Hence, $L$ is split,
which contradicts the hypothesis that it is hyperbolic.

Thus, we have shown that if $M - {\rm int}(N(K))$ contains an essential torus, then there is
one that is invariant under the involution and this descends to a 2-sphere that intersects $L$ in four points.
Since this 2-sphere is separating in $S^3$, we deduce that the essential invariant torus in $M - {\rm int}(N(K))$
is necessarily separating. We pick an essential invariant torus $T$ that is furthest from $K$, in the following sense.
If $T'$ is another essential invariant torus in $M - {\rm int}(N(K))$ that is disjoint from $T$ but not parallel to $T$,
then $T'$ lies in the component of $M - {\rm int}(N(K \cup T))$ containing $\partial N(K)$.
Let $S$ be the image of $T$ in $S^3$. Since $L$ is 2-string prime, this bounds a 3-ball $B$ that contains
a compression disc for $S - L$ disjoint from $L$. Since $T$ was disjoint from $K$, 
its image $S$ is disjoint from $\alpha$. Therefore, $\alpha$ must lie in $B$, because
otherwise the inverse image $T$ of $S$ is compressible in $M - {\rm int}(N(K))$. 
Thus, the tangle $A = B \cap L$ becomes a new 1-manifold $A_\circ = B \cap L_\circ$.

We claim that $\partial B - \partial A_\circ$ is compressible in the complement of $A_\circ$. For if it is incompressible,
then one may find a splitting sphere for $L_\circ$ that is disjoint from it. This splitting sphere cannot lie
in $B$, by Theorem 2.2 (i). Hence, the splitting sphere lies in the complement of $B$, and
therefore forms a splitting 2-sphere for $L$, contrary to assumption. 

Hence, by Theorem 2.2 (ii), $A_\circ$ is
a trivial 2-string tangle. So conclusion (i) or (ii) of Theorem 2.1 holds. We can view the removal
of $A$ and the insertion of $A_\circ$ as tangle replacement along the core arc $\beta$ of $A$.
In conclusion (i) of Theorem 2.1, $\beta$ equals $\alpha$. But in conclusion (ii), $\beta$ is
different from $\alpha$. Let $K'$ be the inverse image of $\beta$ in $M$. Then $M_\circ$
is obtained from $M$ by Dehn surgery along $K'$. The distance between the surgery slope
and the meridian slope is equal to the distance between the slopes of $A$ and $A_\circ$,
and by Theorem 2.1, this is at least $d \geq 2$. Note that by our choice of $T$, $M - {\rm int}(N(K'))$
is atoroidal. Hence, we may apply Theorem 5.1 to $K'$ instead of $K$. 

Thus, we have verified the hypotheses from Theorem 5.1. 

\vskip 18pt
\centerline {\caps 8. The mapping class group of a hyperbolic 3-manifold}
\vskip 6pt

Any branched double cover comes equipped with an involution. Therefore, in this
section, we analyse the {\sl mapping class group} of a compact orientable 3-manifold $X$,
by which we mean the group of homeomorphisms of $X$, up to isotopy.
As we have dealt with the Seifert fibred case in Section 6, 
the manifolds $X$ that we will consider will be hyperbolic. It is well known that the mapping class group of such a
manifold is finite and computable. Indeed, we have the following result.

\noindent {\bf Theorem 8.1.} {\sl Let $X$ be an orientable finite-volume hyperbolic 3-manifold. 
Then the mapping class group of $X$ is finite. Moreover, there is an algorithm that takes, as its input, a triangulation $T$
for $X$ and returns the following:
\item{(i)} a finite sequence of Pachner moves taking $T$ to a triangulation $T'$;
\item{(ii)} a finite group of symmetries of $T'$, which forms a realisation of the mapping
class group of $X$.

}

This result is well known, and this is not the place to explain it in detail. It uses essentially
the same methods as the solution to the homeomorphism problem for compact orientable 3-manifolds.
See for example [22] or [46]. A statement of the
computability of the mapping class group of $X$ is given in Theorem 8.3 of [22] for example,
and the proof there gives (i) and (ii) of Theorem 8.1.

It is also worth pointing out that, in the situations where we want to apply Theorem 8.1,
$\partial X$ is non-empty, and in this case, there is a nice algorithm to solve the problems in Theorem 8.1, as
follows.

It was shown by Petronio and Weeks [43] that when $\partial X$ is a non-empty collection of tori,
$X$ admits a hyperbolic structure if and only if
it has an ideal triangulation that admits a `partially flat' solution to the hyperbolic gluing equations.
Thus, one first transforms the given triangulation into an ideal one.
Then one applies all possible 2-3 and 3-2 Pachner moves to this, to create a list of ideal
triangulations for $X$. Then, for each triangulation in this list, one applies all possible 2-3 and 3-2
Pachner moves, and so on. In this way, an ever-increasing list of ideal triangulations for $X$ is
created. It is a theorem of Matveev [29] that any ideal triangulation for $X$ will eventually appear 
in this list. As this is being produced, the algorithm checks whether each ideal triangulation
admits a partially flat solution to the gluing equations. Note that one can decide whether a given
system of algebraic equations and inequalities with integer coefficients admits a real solution and, if it does,
it is possible to find one. This is due to Tarski [49], although there are now more efficient solutions [13].
If the ideal triangulation does admit a partially flat solution to the gluing equations, such a solution can
therefore be found and this is the required hyperbolic
structure. From this, one can compute the Epstein-Penner decomposition, using the algorithm of 
Weeks [53]. The Epstein-Penner decomposition is a way of building $X$ out of hyperbolic ideal polyhedra
via isometries between their faces. It has the key property that the mapping class group of $X$
is precisely the group of combinatorial automorphisms of these polyhedra that respect the face
identifications. One can then easily decompose the ideal polyhedra into a triangulation $T'$
satisfying the requirements of Theorem 8.1.

The computation of the mapping class group given in Theorem 8.1 can be made effective in the following sense.

\noindent {\bf Theorem 8.2.} {\sl Let $X$ be an orientable finite-volume hyperbolic 3-manifold. Let $h_1$ and $h_2$
be finite order piecewise-linear homeomorphisms of $X$, given as combinatorial automorphisms
of triangulations $T_1$ and $T_2$ of $X$, together with a finite sequence of Pachner moves relating $T_1$ and $T_2$. 
Then there is an algorithm to determine whether $h_1$ and $h_2$ are equal in the mapping class group of $X$.}

\noindent {\sl Proof.} By a theorem of Gabai (Theorem 1.2 in [11]), two homeomorphisms
of a finite-volume hyperbolic 3-manifold $X$ are homotopic if and only if they are isotopic.
Thus, determining whether $h_1$
and $h_2$ are equal in the mapping class group is equivalent to determining whether
$h_2 h_1^{-1}$ is homotopic to the identity. This is equivalent to the induced map
on $\pi_1(X)$ being an inner automorphism. This can be determined as follows.
Pick a generating set $\gamma_1, \dots, \gamma_n$ for $\pi_1(X)$. The hyperbolic
structure on $X$, which has been determined using [22] or [46], realises $\gamma_1, \dots, \gamma_n$
as elements $A_1, \dots, A_n$ of ${\rm PSL}(2, {\Bbb C})$. The homomorphism 
induced by $h_2 h_1^{-1}$ sends $A_1, \dots, A_n$ to $A'_1, \dots, A'_n$ in ${\rm PSL}(2, {\Bbb C})$. If some $A'_i$
is conjugate to $A_i$, the conjugating element would have to send the fixed-point set
for $A'_i$ in the sphere at infinity to the fixed-point set of $A_i$. Hence, once we consider
a couple of loxodromic $A_i$ with disjoint fixed-point sets (which can readily be arranged by a minor
adjustment to the generating set), then
there are only finitely many possible elements of ${\rm PSL}(2, {\Bbb C})$ that can conjugate
the ordered set $A'_1, \dots, A'_n$ to $A_1, \dots, A_n$. If none of these lies in $\pi_1(X)$, then
this can be determined and hence, it can be deduced that the homomorphism induced by 
$h_2 h_1^{-1}$ is not an inner automorphism. On the other hand, if it is an inner automorphism,
then an exhaustive search through the possible conjugating elements of $\pi_1(X)$ will eventually establish
that it is indeed an inner automorphism. $\square$

We also note that equality in the mapping class group for such homeomorphisms
is equivalent to something rather stronger.

\noindent {\bf Theorem 8.3.} {\sl Let $X$ be an orientable finite-volume hyperbolic 3-manifold.
Then two finite order piecewise-linear homeomorphisms $h_1$ and $h_2$ of $X$ are
isotopic if and only if there is a piecewise-linear homeomorphism $\phi$ of $X$ that is isotopic to the identity
and that satisfies $h_2 = \phi^{-1} h_1 \phi$. Moreover, there is an algorithm to find
such a homeomorphism $\phi$ when one is given $h_1$ and $h_2$
as combinatorial automorphisms of triangulations $T_1$ and $T_2$ of $X$,
together with a finite sequence of Pachner moves relating $T_1$ and $T_2$.}

\noindent {\sl Proof.} Let $O_1$ and $O_2$ be the orbifolds $X / \langle h_1 \rangle$ and
$X / \langle h_2 \rangle$. These are orbifold-irreducible and orbifold-atoroidal
since $X$ is irreducible and atoroidal. Hence, they admit hyperbolic structures by the Orbifold Theorem [3, 7].
These lift to hyperbolic metrics $g_1$ and $g_2$ on $X$. By Gabai's theorem
(which is the analogue of the Smale Conjecture for hyperbolic 3-manifolds, Theorem 7.3 [11]),
there is an isotopy between $g_1$ and $g_2$. In other words, there is a 1-parameter
family of diffeomorphisms $\phi_t$ ($t \in [0,1]$) such that $\phi_0$ is the identity and $\phi_1^\ast g_1 = g_2$.
Let $\phi = \phi_1$. Now, $h_1$ and $h_2$ are isometries with respect to $g_1$ and $g_2$ respectively.
Hence, $\phi^{-1} h_1 \phi$ is an isometry with respect to $\phi^\ast g_1 = g_2$.
It is homotopic to the isometry $h_2$ of $g_2$. By Mostow rigidity, two homotopic
isometries are equal, and hence $\phi^{-1} h_1 \phi = h_2$, as required. 

Suppose now that we are given isotopic, finite order, piecewise-linear homeomorphisms $h_1$ and $h_2$ of $X$,
given as combinatorial automorphisms of triangulations $T_1$ and $T_2$ of $X$,
together with a finite sequence of Pachner moves relating $T_1$ and $T_2$.
We then know that $\phi$, satisfying the above conditions, exists. We must give
an algorithm to find it.
Note that $\phi$ descends to a homeomorphism $\overline \phi \colon O_2 \rightarrow O_1$
that respects the singular locus. This satisfies $(\overline \phi)_\ast (p_2)_\ast \pi_1(X) = (p_1)_\ast \pi_1(X)$
where $p_1 \colon X \rightarrow O_1$ and $p_2 \colon X \rightarrow O_2$ are
the quotient maps and $(p_i)_\ast \colon \pi_1(X) \rightarrow \pi_1(O_i)$ are the induced homomorphisms
at the level of orbifold fundamental groups. Conversely, given a homeomorphism $\overline \phi$ between $O_2$ and $O_1$ respecting the singular locus
and satisfying $(\overline \phi)_\ast (p_2)_\ast \pi_1(X) = (p_1)_\ast \pi_1(X)$, we may lift it to a
homeomorphism $\phi \colon X \rightarrow X$ such that $\phi^{-1} h_1 \phi = h_2^k$ for some $k \in {\Bbb Z}$.
Moreover, when $\phi$ is isotopic to the identity, then we may take $k = 1$.
When $h_1$ and $h_2$ are given
as combinatorial automorphisms of triangulations $T_1$ and $T_2$ of $X$,
we may subdivide these to triangulations $T_1'$ (respectively, $T_2'$) with the property
that a simplex is invariant under $h_1$ (respectively $h_2$) if and only if it is in the fixed-point set of $h_1$ 
(respectively $h_2$).  These descend to triangulations of $O_1$ and $O_2$, and we may then start to
search for homeomorphisms $\overline \phi$ between them, by exhaustively trying sequences
of Pachner moves. We will eventually find such a homeomorphism that lifts to a homeomorphism
$\phi \colon X \rightarrow X$ that is isotopic to the identity and that satisfies $\phi^{-1} h_1 \phi = h_2$. $\square$


\vskip 18pt
\centerline {\caps 9. The algorithm to detect links with splitting number one}
\vskip 6pt

In this section, we provide the algorithm required by Theorems 1.1 and 1.2. 

We have already seen that when $L$ admits a tangle replacement along an arc $\alpha$
that creates a split link, then there is an associated exceptional surgery curve $K$
in the manifold $M$. When the distance of the tangle replacement is at least $2$,
this satisfies the hypotheses of Theorem 5.1, and so
this theorem provides a list of all possibilities for $K$ {\sl up to isotopy of $M$}.
But it may not be clear whether a knot $K'$ provided by Theorem 5.1 is isotopic to a curve that
is invariant under the involution on $M$. Even if it is, it is not clear whether
it has several different representatives in its isotopy class in $M$, each of
which is invariant under the involution, but which descend to non-isotopic arcs in 
the exterior of $L$. Thus, it is not immediately clear how to create a finite list of all possibilities
for the arc $\alpha$. To circumvent this problem, we argue as follows. 

The construction of $M$ as a branched double cover provides a piecewise-linear involution $\tau$ 
of $M$ that restricts to a piecewise-linear involution of $M - {\rm int}(N(K))$.
Hence, if $K'$ is isotopic to $K$, then $M - {\rm int}(N(K'))$ also admits a piecewise-linear involution.
So, for each of the knots $K'$ provided by Theorem 5.1, we check whether
$M - {\rm int}(N(K'))$ admits a piecewise-linear involution. The main case that we will consider
is where $M$ and $M - {\rm int}(N(K))$ are hyperbolic. In this situation, their
mapping class groups are finite and computable using Theorem 8.1. So, we
can decide whether $M - {\rm int}(N(K'))$ admits a piecewise-linear involution, and we can find an explicit
representative for each such involution. We are only interested in involutions
that extend to an involution of $M$ that is isotopic to $\tau$. The following
lemma asserts that, if we know the involution of $M - {\rm int}(N(K'))$
and if we know that it extends to the involution $\tau$ on $M$ (up to isotopy), then this is enough to
be able to recreate the image of $K$ in the orbifold $M / \tau$. Recall
that $M / \tau$ is the orbifold with underlying manifold $S^3 - {\rm int}(N(L - L'))$
and with singular set equal to $L'$. The image of $K$ in $M / \tau$ is the required arc $\alpha$.

\noindent {\bf Lemma 9.1.} {\sl Let $\tau$ be a piecewise-linear involution of a hyperbolic 3-manifold $M$
that leaves a knot $K$ invariant, where $M - {\rm int}(N(K))$ is hyperbolic. Let $\rho$ be a piecewise-linear homeomorphism of 
$M$ that is isotopic to the identity, taking $K$ to a knot $K'$. Let $\eta$ be another piecewise-linear involution of $M$ that leaves
$K'$ invariant. Suppose that the restrictions of $\rho \tau \rho^{-1}$ and
$\eta$ are isotopic homeomorphisms of $M - {\rm int}(N(K'))$. Suppose also that 
there is a piecewise-linear homeomorphism $\phi$ of $M$ that is isotopic to the identity and that satisfies $\eta = \phi^{-1} \tau \phi$. 
Then there is a piecewise-linear homeomorphism $M / \tau \rightarrow M / \tau$ respecting the singular locus of this orbifold and taking
$K / \tau$ to $\phi(K') / \tau$. Moreover, this homeomorphism is isotopic to the identity
on the components of $\partial M / \tau$ that are disjoint from the singular set.}

\noindent {\sl Proof.} Since $\rho \tau \rho^{-1}$ and
$\eta$ are isotopic piecewise-linear homeomorphisms of $M - {\rm int}(N(K'))$, then by Theorem 8.3, there is
a piecewise-linear homeomorphism $\psi$ of $M - {\rm int}(N(K'))$, isotopic to the identity, such that 
$\rho \tau \rho^{-1} = \psi^{-1} \eta \psi$. This extends to a piecewise-linear homeomorphism
$\psi$ of $M$ such that $\psi(K') = K'$. Thus, we have the following commutative diagram:

$$

\matrix{
(M,K) & \buildrel \rho \over \longrightarrow & (M,K') & \buildrel \psi \over \longrightarrow & (M, K') & \buildrel \phi \over \longrightarrow & (M, \phi(K'))\cr
\Big\downarrow \rlap{$\vcenter{\hbox{${\scriptstyle \tau}$}}$}
&& \Big\downarrow \rlap{$\vcenter{\hbox{${\scriptstyle \rho \tau \rho^{-1}}$}}$}
&& \Big\downarrow  \rlap{$\vcenter{\hbox{${\scriptstyle \eta}$}}$}
&& \Big\downarrow \rlap{$\vcenter{\hbox{${\scriptstyle \tau}$}}$} \cr
(M,K) & \buildrel \rho \over \longrightarrow & (M,K') & \buildrel \psi \over \longrightarrow & (M, K') & \buildrel \phi \over \longrightarrow & (M, \phi(K'))\cr
}
$$
Hence, the composition $\phi \psi \rho$ descends to a homeomorphism
$M / \tau \rightarrow M / \tau$ taking singular locus to singular locus and taking
$K / \tau$ to $\phi(K') / \tau$. Since $\phi \psi \rho$ is isotopic to the identity, its action on each component of $\partial M / \tau$ that is disjoint from the singular set
is isotopic to the identity. $\square$

Thus, the algorithm required by Theorems 1.1 and 1.2 is as follows:

\item{1.} If $L$ is provided by a diagram, then we use this to build a triangulation
of $S^3$ in which $L$ is simplicial.
\item{2.} Pick all sublinks $L'$ of $L$ consisting of one or two components. 
If $|L| = 2$, then we require that $|L'| = 1$. There are only finitely many
choices for $L'$ and so let us focus on just one such choice.
\item{3.} Construct a triangulation $T$ of the double cover of $S^3 - {\rm int}(N(L - L'))$ branched over $L'$.
Denote this manifold by $M$. Note that $M$ has non-empty boundary.
\item{4.} By Proposition 4.1, $M$ is either hyperbolic or a small Seifert fibre space.
Using Algorithm 8.1 in [19], determine which of these cases holds. The algorithm divides into these two cases.

\noindent {\caps The hyperbolic case}

\item{1.} Use Theorem 5.1 to produce a finite list of knots $K$ in $M$, with slopes
$\sigma$ on $\partial N(K)$ that satisfy the hypotheses of Theorem 5.1.
Each knot $K$ can be given as a subcomplex of some iterated barycentric subdivision
of $T$. Let us now focus on just one choice of $K$ and $\sigma$.
\item{2.} From the way that $K$ is given, it is easy to build a triangulation
$T'$ for $M - {\rm int}(N(K))$. For example, one can take two further barycentric subdivisions
of the triangulation of $M$, and then remove the simplices that are incident to $K$.
\item{3.} Use Theorem 4.1.12 in [30] or Theorem 5.2 in [19] to determine whether $M - {\rm int}(N(K))$ is irreducible
and use Theorem 6.4.10 in [30] or Algorithms 8.1 and 8.2 in [19] to determine whether it is atoroidal.
If it is reducible or toroidal, then discard it and move on to the next choice
of $K$ and $\sigma$. So let us assume that $M - {\rm int}(N(K))$ is irreducible and atoroidal. 
It is therefore hyperbolic or Seifert fibred, by the solution to the Geometrisation Conjecture.
In fact, it cannot be Seifert fibred, because it can be Dehn filled to form the hyperbolic manifold $M$.
\item{4.} Compute the mapping class group of $M - {\rm int}(N(K))$ using Theorem 8.1.
This gives a triangulation $T''$ of $M - {\rm int}(N(K))$ and a group of symmetries
of $T''$ that realises the mapping class group. It also provides a sequence of Pachner moves
from $T'$ to $T''$.
\item{5.} For each order 2 symmetry $\eta$ of $M - {\rm int}(N(K))$, determine
whether it preserves $\partial N(K)$ and determine whether it
acts by $-{\rm id}$ on it. If this is not the case, ignore it and move on.
\item{6.} If $\eta$ does act in this way on $\partial N(K)$, then it
extends to an order two symmetry of $M$, which we will also call $\eta$. 
Extend $T''$ to a triangulation $T'''$ of $M$ that is invariant under $\eta$ and build a sequence
of Pachner moves from $T$ to $T'''$.
\item{7.} Using Theorem 8.2, determine whether $\eta$
is isotopic to the involution $\tau$ of $M$ that is given by the construction of $M$
as a branched double cover. If it is not, then discard it.
\item{8.} Assuming that $\eta$ is isotopic to $\tau$,
Theorem 8.3 provides a piecewise-linear homeomorphism $\phi$ of $M$ isotopic to the identity
such that $\eta = \phi^{-1} \tau \phi$. The arc $\phi(K)$
is invariant under $\tau$ and its image in $S^3 - {\rm int}(N(L - L'))$
is an arc $\beta$ with endpoints on $L'$.
\item{9.} Lemma 9.1 only provides the arc $\beta$ up to homeomorphism of $S^3 - {\rm int}(N(L))$, 
whereas we want all possibilities for $\beta$ up to isotopy.
So, for each arc $\beta$ constructed as above, consider all its images under the mapping class
group of $S^3 - {\rm int}(N(L))$, which is determined using Theorem 8.1.
\item{10.} Feed all these arcs $\beta$ into the subroutine below that constructs
arcs $\alpha$ from them.

\noindent {\caps The Seifert fibred case}

\item{1.} Use Theorem 6.6 to produce a Seifert fibration of $M$ that is invariant under the involution.
More specifically, this produces a union of disjoint simple closed curves $K$ in $M$ that are the singular fibres.
Also, for each such simple closed curve, it produces the slope $\sigma$ on $\partial N(K)$ that is slope
of the regular fibres.
\item{2.} For each possibility for $K$, observe whether it is invariant under the involution. If it is not, discard it.
If it is, its image in $S^3 - {\rm int}(N(L - L'))$
is an arc $\beta$. Feed $\beta$ into the subroutine below that constructs
arcs $\alpha$ from it.

\noindent {\caps Constructing the arcs $\alpha$}

\item{1.} Using $\beta$, we determine one or two possibilities for $\alpha$,
corresponding to conclusions (i) and (ii) of Theorem 2.1. In case (i),
we set $\alpha$ to be equal to $\beta$. We check that the distance of the tangle
replacement is two. If it is not, we discard this possibility. In case (ii),
we parametrise the slopes on $\partial N(K)$ by ${\Bbb Q} \cup \{ \infty \}$
where $\infty$ is the meridional slope giving $M$. We let the slope
of $\sigma$ be $p/q$. If $q$ is of the form $2 a^2$, then we write 
$p/q$ as $(1 \pm 2ab) / \pm 2a^2$, and then set $\alpha$ as in (ii) of Theorem 2.1.
\item{2.} If we are determining whether $s_d(L) = 1$ or $ts_d(L)=1$,
then we only consider arcs $\alpha$ that have endpoints on distinct components
of $L$. So under these circumstances, if the endpoints of $\alpha$ lie
on the same component of $L$, then we discard it.
\item{3.} For each of these possibilities for $\alpha$, we perform this tangle replacement. The image of $L$ is a link $L_\circ$.
Determine whether $L_\circ$ is the unlink, a split link or totally split, as appropriate.
For example, one can use Theorem 5.2 in [19] to determine whether $S^3 - {\rm int}(N(L_\circ))$ is reducible
and if it is, to find a reducing sphere. Then one can decompose along it, fill in with 3-balls
and repeat. In this way, we can determine whether $L_\circ$ is split and whether it is
totally split. In the latter case, we can also determine whether the components are
unknots, and so whether $L_\circ$ is the unlink.
\item{4.} For any relevant crossing arc $\alpha$, we can then construct the associated
crossing circles and $\pm 1$ surgery coefficients, as in Section 3.

One may wonder why we could not simply set $\alpha$ to be $\beta$ in Step 1 of the final subroutine.
Note that $\alpha$ is the image of $K$ under the quotient map $M \rightarrow M / \tau$.
At the end of Section 7, it was necessary to permit $K$ to be replaced by another knot $K'$.
In this case, $\beta$ is the image of $K'$, and $\alpha$ is obtained from it by the procedure
described in Theorem 2.1.

Note that using a minor variation of Step 1 in the final subroutine, we can generalise from crossing
changes to tangle replacements with distance at least two. The modified version of
Step 1 is as follows. In case (i) of Theorem 2.1, we set $\alpha$ equal to $\beta$, but
we do not check that the distance of the tangle replacement is two. This was done solely
to ensure that the tangle replacement corresponded to a crossing change. In case (ii)
of Theorem 2.1, we consider all possible ways of writing $p/q$ as $(1 \pm dab) / \pm da^2$,
where $a$ and $b$ are integers and $d$ is an integer at least two. For each such possibility,
we get an arc $\alpha$ as in (ii) of Theorem 2.1. Thus, checking each of these
tangle replacements in turn, we obtain the following result.

\noindent {\bf Theorem 9.2.} {\sl Let $L$ be a link in $S^3$ with at least two components.
If $L$ has exactly two components, suppose that these have zero linking number.
Suppose that $L$ is hyperbolic and 2-string prime. Then there is an algorithm to find all possible
trivial tangle replacements that can be made to $L$ with distance at least two that turn
it into a split link.

}

\noindent {\bf Remark 9.3.} The algorithm given above required us to consider all sublinks $L'$
of $L$ with one or two components. However, in the case where $|L| \geq 3$, we can focus just
on the case where $|L'| = 2$. The reason for this is as follows. The manifold $M$ is the
double cover of $S^3 - {\rm int}(N(L - L'))$ branched over $L'$. It was important for the
endpoints of $\alpha$ to lie in $L'$, so that the inverse image of $\alpha$ in $M$ is a knot $K$.
But if $\alpha$ has endpoints in the same component of $L$, then we can choose $L'$ to
be the union of this component plus one other chosen arbitrarily. The conditions on $M - {\rm int}(N(K))$
are all easily verified, as in Section 7. In particular, the condition $H_2(M - {\rm int}(N(K)),\partial M) \not= 0$
holds.

\vskip 18pt
\centerline{\caps 10. The algorithm to enumerate knots with exceptional surgeries}
\vskip 6pt

In the previous section, an algorithm that solves the decision problems in Theorem 1.1 
was given. The crucial ingredient was Theorem 5.1, which provided an algorithm
to enumerate the exceptional surgery curves within a 3-manifold satisfying certain
conditions. This algorithm is difficult to work with in practice, and so in many concrete
examples, it is better to use the techniques behind Theorem 5.1.
Therefore in this section, we give an overview of the proof of Theorem 5.1.

The proof relied heavily on sutured manifold theory. An excellent reference for this is [44].
The first part of the argument closely follows Section 5 of [44].

Suppose that $M$ is a compact orientable 3-manifold with $\partial M$ a (possibly empty) union of tori.
Let $K$ be a knot in $M$ and let $\sigma$ be a slope on $\partial N(K)$ satisfying (i) and (ii) of Theorem 5.1.
Give $M - {\rm int}(N(K))$ the structure of a sutured manifold with $R_- = \emptyset$ and $R_+ = \partial M \cup \partial N(K)$ and
therefore with sutures $\gamma_1 = \emptyset$.
We may find a taut sutured manifold hierarchy
$$(M - {\rm int}(N(K)), \gamma_1) \buildrel S_1 \over \longrightarrow (M_2 - {\rm int}(N(K)), \gamma_2) \buildrel S_2 \over \longrightarrow
\dots \buildrel S_{n-1} \over \longrightarrow (M_n - {\rm int}(N(K)), \gamma_n)$$
such that the following hold:
\item{(i)} Each surface $S_i$ is disjoint from $\partial N(K)$.
\item{(ii)} Each surface $S_i$ contains no closed separating components.
\item{(iii)} No surface $S_i$ has a boundary component that bounds a disc in $\partial M_i$ disjoint from $\gamma_i$.
\item{(iv)} $H_2(M_n - {\rm int}(N(K)), \partial M_n) = 0$.
\item{(v)} In the case where $\sigma$ is norm-exceptional there is some
$z \in H_2(M - {\rm int}(N(K)), \partial M)$ that maps to an element 
$z_\sigma \in H_2(M_K(\sigma), \partial M_K(\sigma))$, such that the
Thurston norm of $z_\sigma$ is less than the Thurston norm of $z$. We require that $[S_1] = z$.

Condition (iv) implies that $\partial M_n$ consists of a collection of spheres and at most one torus. By the
tautness of $(M_n - {\rm int}(N(K)), \gamma_n)$, each sphere bounds a ball in $M_n$. Hence, there
is exactly one component $Y$ of $M_n - {\rm int}(N(K))$ that is not a ball. This forms a rational homology cobordism
between a toral component $T$ of $\partial M_n$ and $\partial N(K)$. This torus $T$ is incompressible in $M - {\rm int}(N(K))$.
This is because a compressible torus in an irreducible 3-manifold
bounds a solid torus or lies within a 3-ball.
If $T$ bounds a solid torus in $M - {\rm int}(N(K))$, then some surface $S_i$ must intersect this
solid torus, and this gives rise to a closed separating component of $S_i$, contradicting (ii). If
$T$ lies within a 3-ball in $M - {\rm int}(N(K))$, then again some $S_i$ must intersect this 3-ball
and again this gives a closed separating component. Now we are assuming that $M - {\rm int}(N(K))$
is atoroidal, and so the incompressible torus $T$ is boundary parallel in $M - {\rm int}(N(K))$. It cannot be parallel to
component of $\partial M$, because $Y$ would then be the region between $T$ and $\partial N(K)$
and hence would be a copy of $M - {\rm int}(N(K))$. However, $M - {\rm int}(N(K))$ is not
a rational homology cobordism between $\partial M$ and $\partial N(K)$, since 
$H_2(M - {\rm int}(N(K)), \partial M)$ is non-trivial. We therefore deduce that
$T$ is parallel to $\partial N(K)$, and hence $Y$ is a copy of $T \times [0,1]$ with
$T \times \{ 1 \} = T$ and $T \times \{ 0 \} = \partial N(K)$. 

Consider any slope $\rho$ on $\partial N(K)$ other than the one that is parallel to
the sutures $\gamma_n \cap Y$. If we Dehn fill each of the
manifolds $M_i - {\rm int}(N(K))$ along this slope, we obtain a sutured manifold hierarchy
$$(M(\rho), \gamma_1) \buildrel S_1 \over \longrightarrow (M_2(\rho), \gamma_2) \buildrel S_2 \over \longrightarrow
\dots \buildrel S_{n-1} \over \longrightarrow (M_n(\rho), \gamma_n).$$
Here, $M_i(\rho)$ denotes the result of performing Dehn surgery along $K$ in $M_i$ with slope $\rho$.
Since we assumed that $\rho$ is not parallel to
the sutures $\gamma_n \cap Y$, the sutured solid torus $(M_n(\rho), \gamma_n)$ is taut.
Using the theorem that tautness pulls back (Theorem 3.6 of [44]), we deduce that each of the manifolds
in the hierarchy is taut and each of the decomposing surfaces is taut. Because $\sigma$ is exceptional
or norm-exceptional, $(M(\sigma), \gamma_1)$ is not taut or $S_1$ is not taut in $M(\sigma)$.
Hence, $\sigma$ must be the slope on $\partial N(K)$ parallel to the sutures $Y \cap \gamma_n$.
So, if we re-attach the solid torus $N(K)$ using the meridional slope, we deduce that
$$(M, \gamma_1) \buildrel S_1 \over \longrightarrow (M_2, \gamma_2) \buildrel S_2 \over \longrightarrow \dots
\buildrel S_{n-1} \over \longrightarrow (M_n, \gamma_n)$$
is taut.

The key part of the proof of Theorem 5.1 is to place this hierarchy into some sort of `normal form'
with respect to a given triangulation ${\cal T}$ of $M$. In fact, we first dualise ${\cal T}$ to form a handle structure
${\cal H}$ for $M$, and we make the hierarchy `normal' with respect to ${\cal H}$.
We will shortly make this statement a little more precise,
but the idea is roughly that there should be only finitely possibilities for $M_i \cap H$ and $\gamma_i \cap H$,
for each handle $H$ of ${\cal H}$. In fact, this statement is not quite correct, but we will
make it accurate shortly. But if there were only finitely many possibilities for $M_n \cap H$ and $\gamma_n \cap H$,
and we could enumerate them, then we could reconstruct the way that $(M_n, \gamma_n)$ lies within
$M$. Since $M_n$ is a regular neighbourhood of $N(K)$, this would imply that there are only
finitely many possibilities for $K$ and $\sigma$, and we could enumerate them all.

To make the above discussion more precise, it is helpful to consider just the first decomposition
along $S_1$. Since $S_1$ is taut in $M$, it is incompressible. Suppose also that $S_1$ is also boundary-incompressible. Therefore,
we may place $S_1$ into normal form with respect to the handle structure ${\cal H}$. 
(A notion of normal surfaces in a handle structure was first defined by Haken [14]; see also
[17] and [30]. A variant of this notion was in fact used in [24], partly to take account of the possibility that $S_1$ may be boundary-compressible.)
Then when we cut along $S_1$, we obtain a handle structure ${\cal H}_2$ for $M_2$.
It would be convenient if, within each handle $H$ of ${\cal H}$, there are only finitely
many possibilities for $M_2 \cap H$. However, this need not be the case. The surface
$S_1$ may have many components of intersection with $H$, and thereby give rise to
many handles of ${\cal H}_2$ within $H$. To get around this problem, we use the notion
of the parallelity bundle of ${\cal H}_2$. By definition, a handle of ${\cal H}_2$ is a {\sl parallelity
handle} if it lies between two normally parallel discs of $S_1$. This is an $I$-bundle, with
$\partial I$-bundle lying in the copies of $S_1$ in $\partial M_2$. The union of the parallelity
handles is an $I$-bundle ${\cal B}$ over a surface $F$ called the {\sl parallelity bundle} for ${\cal H}_2$. 
It is clear that, within each
handle $H$ of ${\cal H}$, there are only finitely many possibilities for $H \cap (M_2 - {\cal B})$ and
these are algorithmically constructible. This list is universal, in the sense that it does not depend
on $M$ or any other data. It only depends on the way that $H$ intersects the neighbouring handles
of higher index, and there are only finitely many possibilities for this because ${\cal H}$ is dual to
a triangulation.

One of the key parts of the proof of Theorem 5.1 is therefore to remove the parallelity bundle
${\cal B}$. The procedure is given in detail in Section 8 of [24]. It involves making changes
to the handle structure ${\cal H}_2$. We now explain the most important of these changes now.
The {\sl horizontal boundary} $\partial_h {\cal B}$ is the $(\partial I)$-bundle
and lies in the copies of $S_1$. The {\sl vertical boundary} $\partial_v {\cal B}$ is the $I$-bundle over $\partial F$.
The surface $P = {\rm cl}(\partial_v {\cal B} - \partial M)$ is also an $I$-bundle, and hence
it is a collection of discs and annuli. It is a properly embedded surface in $M_2$.
A key modification that is made is to decompose $M_2$ along some of the components
of $P$. Each such component is either an annulus disjoint from the sutures $\gamma_2$ or a {\sl product disc},
which is a disc intersecting $\gamma_2$ twice. This surface is
properly embedded in $M_2$, but there is no {\sl a priori} reason why it should be disjoint from $K$.
This is a consequence of the following result, which was
the central result of [23]. It was here that the hypothesis that the distance $\Delta(\mu, \sigma) > 1$ is used.

\noindent {\bf Theorem 10.1.} {\sl Let $(M, \gamma)$ be a taut sutured manifold. Let
$K$ be a knot in $M$ such that $M - {\rm int}(N(K))$ is irreducible and atoroidal. Let 
$\sigma$ be a slope on $\partial N(K)$ such that $\Delta(\sigma, \mu) > 1$,
where $\mu$ is the meridional slope. Let $M_K(\sigma)$ be the result of
performing Dehn surgery along $K$ with slope $\sigma$. Suppose that
$(M_K(\sigma), \gamma)$ is not taut. Let $G$ be a surface properly embedded in $M$ with components $G_1, \dots, G_{|G|}$,
none of which is a sphere or disc disjoint from $\gamma$. Then 
there is an ambient isotopy of $K$ in $M$ after which, for each integer $i$ between $1$ and $|G|$, we have
$$|K \cap G_i| \leq {-2 \chi(G_i) + |G_i \cap \gamma| \over 2 (\Delta(\mu, \sigma) - 1)}.$$
}

Note that when $G$ is a union of annuli disjoint from $\gamma$ and product discs, then Theorem 10.1 implies
that the knot $K$ may be ambient isotoped off $F$. The condition that $M - {\rm int}(N(K))$ is atoroidal
in fact can be weakened somewhat (see Theorem 1.4 in [23]).

Thus, the proof of Theorem 5.1 proceeds as follows. We start with a taut decomposition
$$(M - {\rm int}(N(K)), \gamma_1) \buildrel S_1 \over \longrightarrow (M_2 - {\rm int}(N(K)), \gamma_2)$$
satisfying (i)-(v) above.
As argued above, when we attach the solid torus using the meridional Dehn filling, we get
a taut decomposition
$$(M, \gamma_1) \buildrel S_1 \over \longrightarrow (M_2, \gamma_2).$$
Because it is taut, $S_1$ can be placed in a position rather similar to normal form with respect to ${\cal H}$ (specifically,
it satisfies Conditions 1-5 of Section 9 in [24]). In fact, one may need to modify $S_1$ to place it in
this form; we will discuss this below. Let ${\cal B}$ be the parallelity bundle in $M_2$ associated with
$S_1$. As discussed above, the surgery curve $K$ can be isotoped off its vertical boundary,
and hence off ${\cal B}$ altogether. We can then apply the procedure given in Section 8 of [24]
to modify ${\cal B}$, to give a new 3-manifold $(M_2', \gamma_2')$. This modification has the
effect of removing all components of ${\cal B}$ that are not $I$-bundles over discs and all components
for which the interior of the vertical boundary intersects $\partial M_2$. The remaining components of 
${\cal B}$, which are therefore $I$-bundles over discs, become 2-handles of $(M'_2, \gamma_2')$.
We have little control over the location of these 2-handles, but we do have control over their
attaching locus onto the 0-handles and 1-handles. This manifold $(M'_2, \gamma_2')$ contains $K$,
and it has the property that when surgery along $K$ is performed, the resulting sutured manifold
is not taut. The advantage of working with this new manifold $M'_2$ is that, for each
handle $H$ of ${\cal H}$, there are only finitely many possibilities for the intersection between
$H$ and the 0-handles and 1-handles of $M'_2$, the attaching locus of the 2-handles and the sutures $\gamma_2'$.

As mentioned above, before $S_1$ satisfies Conditions 1-5 of Section 9 in [24], it may be necessary
to make some modifications to it. These are given in Section 9 of [24]. At each stage, it is ensured that
$K$ remains disjoint from $S_1$. For example, $S_1$ may be boundary-compressed along a product disc.
Using Theorem 10.1, we can ensure that $K$ avoids this product disc and so remains disjoint from
the new surface.

One can then repeat this procedure. We find a taut decomposition 
$$(M'_2 - {\rm int}(N(K)), \gamma_2') \buildrel S_2 \over \longrightarrow (M_3 - {\rm int}(N(K)), \gamma_3)$$
satisying (i)-(iv) above.
Since surgery along $K$ with slope $\sigma$ gives a sutured manifold that is not taut, when we
Dehn fill along the meridional slope, we get a taut decomposition
$$(M'_2, \gamma_2') \buildrel S_2 \over \longrightarrow (M_3, \gamma_3).$$
Now isotope $K$ off the parallelity bundle in $(M_3, \gamma_3)$, and
then remove this bundle to get a new sutured manifold $(M'_3, \gamma'_3)$
containing $K$.

We end with a manifold $(M'_n, \gamma'_n)$ containing $K$ such that
$H_2(M'_n - {\rm int}(N(K)), \partial M'_n)$ is trivial. Hence, as argued above,
$M'_n$ is some 3-balls plus a solid torus with $K$ as its core curve, and the exceptional slope
$\sigma$  on $\partial N(K)$ is parallel to the sutures on $\partial M'_n$.
For each handle $H$ of ${\cal H}$, its intersection with the 0-handles of $M'_n$, the
1-handles of $M'_n$, the attaching locus of the 2-handles and the sutures $\gamma'_n$ takes
one of only finitely many possibilities. The algorithm proceeds by inserting all such possibilities
into each handle of ${\cal H}$ so that they patch together correctly along adjacent handles.
In this way, we can build all possibilities for the sutured manifold $(M'_n, \gamma'_n)$.
The algorithm checks, for each possible $(M'_n, \gamma'_n)$ whether it is a taut solid
torus plus possibly some taut 3-balls. If it is, the core curve of the solid torus is
a possibility for $K$ and the slope of the sutures on the solid toral component of 
$(M'_n, \gamma'_n)$ is a possibility for the slope $\sigma$. Once one has this list,
one can determine easily whether $M$ and $K$ really do satisfy (i) and (ii) of Theorem 5.1.

We summarise the above discussion in the following theorem.

\vfill\eject
\noindent {\bf Theorem 10.2.} {\sl There is a finite computable list of 4-tuples $(H_i, {\cal F}^0_i, {\cal F}^1_i, \gamma_i)$
where
\item{(i)} each $H_i$ is a collection of balls embedded within a tetrahedron $\Delta$;
\item{(ii)} each ${\cal F}^0_i$ is the intersection between $H_i$ and $\partial \Delta$; it is a collection of discs lying in the interior of the faces of $\Delta$;
\item{(iii)} each ${\cal F}^1_i$ is a collection of disjoint rectangles lying in $\partial H_i$; two opposite sides of each rectangle lie in $\partial {\cal F}^0$ and the remainder
of the rectangle is disjoint from ${\cal F}^0$;
\item{(iv)} each $\gamma_i$ is a collection of disjoint arcs properly embedded in ${\rm cl}(\partial H_i - ({\cal F}^0_i \cup {\cal F}^1_i))$.

\noindent These have the following property. Suppose that $M$ is a compact orientable 3-manifold with boundary a (possibly empty) union of tori,
and that $K$ is knot in $M$ with an exceptional or norm-exceptional slope $\sigma$ on $\partial N(K)$, satisfying the hypotheses of
Theorem 5.1. Then for any triangulation of $M$, one may form a handle structure on $N(K)$ as follows. Its 0-handles are
obtained by inserting some $H_i$ into each tetrahedron of the triangulation. The 1-handles are dual to the discs ${\cal F}^0_i$. The
rectangles ${\cal F}^1_i$ patch up to form annuli, which are the attaching locus of the 2-handles. The arcs $\gamma_i$ patch
together to form curves with slope $\sigma$ on $\partial N(K)$.

}

The finite list of 4-tuples in the above theorem is universal, in the sense that it does not depend
on $M$ or $K$. An algorithm to construct this list is given in Section 11 of [24]. 

Thus, these 4-tuples patch together to form a handle structure on $N(K)$. If one wanted to, one
could then realise $K$ in $M$ by picking a curve on $\partial N(K)$ with winding number $1$ in $N(K)$.
This could then be realised as a simplicial curve in some iterated barycentric subdivision of the triangulation of $M$.

The above techniques also provide a proof of Theorem 6.1, which gives that certain exceptional
surgery curves in a Seifert fibred space must be isotopic to an exceptional fibre.

\noindent {\sl Proof of Theorem 6.1.} Let $M$ be a Seifert fibre space with non-empty boundary.
Let $K$ be a knot in $M$ and let $\sigma$ be a slope on $\partial N(K)$ satisfying the hypotheses Theorem 6.1.
Give $M - {\rm int}(N(K))$ the structure of a sutured manifold with $R_- = \emptyset$ and $R_+ = \partial M \cup \partial N(K)$ and
therefore with sutures $\gamma_1 = \emptyset$.

We are assuming that $M$ has non-empty boundary.  Therefore, there is a (possibly empty) union of disjoint properly embedded arcs
in its base space that avoid the exceptional points and that decompose the base space either into a collection of
regular neighbourhoods of the exceptional points or, in the case where there are no singular points, into a single
disc. The inverse image of these arcs is a union of disjoint properly embedded annuli $A$ in $M$, such that $M - {\rm int}(N(A))$ is either a regular neighbourhood of
the singular fibres or, in the case where $M$ has no singular fibres, a fibred solid torus. By Theorem 10.1,
there is an ambient isotopy taking $K$ off $A$. Hence, $K$ lies in $M - {\rm int}(N(A))$. By the irreducibility and atoroidality
of $M - {\rm int}(N(K))$, $K$ must be a core curve of one of the components of $M - {\rm int}(N(A))$. Hence,
$K$ is isotopic to a fibre of $M$ and therefore $M - {\rm int}(N(K))$ is Seifert fibred. The exceptional slope
$\sigma$ on $\partial N(K)$ must be the slope of the regular fibres.
Since we are assuming that the distance between $K$ and the meridian is more than $1$, we deduce
that $K$ must be an exceptional fibre of $M$. $\square$

\vskip 18pt
\centerline{\caps 11. Finiteness of the number of splitting crossing changes}
\vskip 6pt

In this section, we prove Theorem 1.2.

\noindent {\sl Proof.} The algorithmic part of Theorem 1.2 was dealt with in Section 9. So we need
only establish the required upper bound on the number of splitting crossing changes that
can be applied to the given link $L$.

We are given a triangulation $T$ of $S^3$ with $t$ tetrahedra in which $L$
is simplicial. Note that if we
are given, alternatively, a diagram of $L$ with $c$ crossings, then we can easily construct
such a triangulation where $t \leq 24c$. One way of doing this is as follows. 

First apply type $1$ Reidemeister moves to remove any edges in the diagram that
start and end at the same crossing. Then place an octahedron
at each crossing of the diagram. The over-arc and the under-arc at the crossing will be
subcomplexes of these octahedra. Lying above the plane of the diagram and all these octahedra
is a 3-ball. Its boundary has a cell structure. Any 2-cell of this cell structure that is not already a
triangle may be subdivided into triangles. We then triangulate the ball by placing a vertex in its interior
and coning off. We triangulate each octahedron using 4 tetrahedra. The 3-ball lying below the plane of the diagram
and the octahedra is triangulated also by coning off its boundary. The result is a triangulation of the 3-sphere
with $L$ has a subcomplex. It is easy to check that at most $24c$ tetrahedra have been used.

We apply the algorithm given in Section 9 to this triangulation, but skipping some steps that
are not relevant to the counting argument. Step 1 has already been completed.

In Step 2, one considers all sublinks $L'$ of $L$ consisting of one or two components.
If $|L| = 2$, then we require that $|L'| = 1$.
The number of such sublinks is at most $|L|(|L| + 1)/2$. The number of components of $L$
is at most the number of 1-simplices of $T$, which is at most $6t$. So,
the number of relevant sublinks is at most a quadratic function of $t$.

In Step 3, the double cover $M$ of $S^3 - {\rm int}(N(L - L'))$ branched
over $L'$ is constructed. It is straightforward to build a triangulation of $M$, starting
from the triangulation of $S^3$ with $L$ as a subcomplex. The number $t'$ of tetrahedra 
in this triangulation can easily be arranged to be at most a linear function of $t$.

The algorithm now divides into the cases where $M$ is hyperbolic or Seifert fibred.
We consider the hyperbolic case first.

The construction of the knots $K$ in $M$ provided by Theorem 5.1 produces
at most $(k_1)^{t'}$ possibilities for $K$ and $\sigma$, where $k_1$ is a universal
computable constant. Specifically, suppose that Theorem 10.2 provides a list of
$k_1$ 4-tuples. Then each possibility for $N(K)$ is obtained by inserting the 0-handles
in one of the 4-tuples into each tetrahedron of the triangulation, in such a way
they patch together correctly along the faces. Moreover, the arcs in the 4-tuples
patch together to form a representative for $\sigma$. Thus, there are at most $(k_1)^{t'}$
possibilities for $N(K)$ and $\sigma$.

For each possible $K$, one can build a triangulation $T'$ for 
$M - {\rm int}(N(K))$. This is Step 2 of the hyperbolic case. Rather than using iterated
barycentric subdivisions, we do this as follows. We subdivide the 0-handles from Theorem 10.2
into tetrahedra and then extend this triangulation over the 1-handles
and 2-handles of $N(K)$. Thus, the number of tetrahedra is
at most $k_2 t$ for a universal computable $k_2$.

In Step 3 of the hyperbolic case, $M - {\rm int}(N(K)$ is discarded
if it is not irreducible and atoroidal.

In Step 4 of the hyperbolic case, the symmetry group of $M - {\rm int}(N(K))$ is computed.
The size of the symmetry group for a hyperbolic 3-manifold $X$ is
at most a linear function $k_3 {\rm vol}(X)$ for the following reason.
The quotient of $X$ by its symmetry group is a finite-volume hyperbolic
orbifold and there is a universal lower bound $v$ on the volume of such
an orbifold [20]. Thus, the order of the symmetry group of $X$ is at most
${\rm vol}(X) / v$. Setting $k_3 = 1/v$ establishes the claim.

Note that $k_3 {\rm vol}(M - {\rm int}(N(K))) \leq k_3 v_3 k_2 t$, where $v_3$ is the volume of
a regular ideal hyperbolic 3-simplex. This follows from the general result [50]
that the volume of hyperbolic 3-manifold with (possibly empty) toroidal boundary
is at most $v_3$ times the number of tetrahedra in any triangulation of the manifold.
Thus, the symmetry group of $M - {\rm int}(N(K))$ has order at most
$k_3 v_3 k_2 t$. Each order two symmetry produces at most one possibility
for the arc $\beta$. 

In Step 9 of the hyperbolic case, we consider all the images of these arcs
$\beta$ under the mapping class group of $S^3 - {\rm int}(N(L))$. This mapping class
group has order at most $k_3 {\rm vol}(S^3 - {\rm int}(N(L))) \leq k_3 v_3 t$.

We now consider the case where $M$ is Seifert fibred. By Theorem 6.1, there is
a Seifert fibration of $M$ in which $K$ is a singular fibre and by Addendum 6.2, this
Seifert fibration can be chosen to be invariant under $\tau$.
Since $M$ is atoroidal and has non-empty boundary, it has at most
two singular fibres. Thus, in the Seifert fibred case, there are
at most two possibilities for the image $\beta$ of $K$.

Thus, in both the case where $M$ is hyperbolic and where it is Seifert fibred,
we have a bound on the number of possibilities for the arc $\beta$ and the associated
tangle replacement slope. For each $\beta$ and associated slope, there are at most two possible arcs $\alpha$.
For each tangle replacement along $\alpha$, there are two associated crossing circles,
by Lemma 3.1. The number of possible crossing circles is therefore at most
$$12t(6t + 1) k_1^{t'} k_3^2 v_3^2 k_2 t^2$$
which is at most $k^t$ for some universal computable constant $k$. 
$\square$



\vskip 18pt
\centerline {\caps 12. The Whitehead link}
\vskip 6pt

In this section, we examine an example, the Whitehead link.
We determine the complete set of crossing changes that turn the link into a
split link. 

\vfill\eject
\noindent {\bf Theorem 12.1.} {\sl Any crossing change that turns the Whitehead link
into a split link is equivalent to changing one of the specified crossings in Figure 8.
In particular, there are $3$ crossing circles up to equivalence, and $2$ crossing arcs
up to equivalence, that yield splitting crossing changes.}

\vskip 12pt
\centerline{
\epsfxsize = 2.5in
\epsfbox{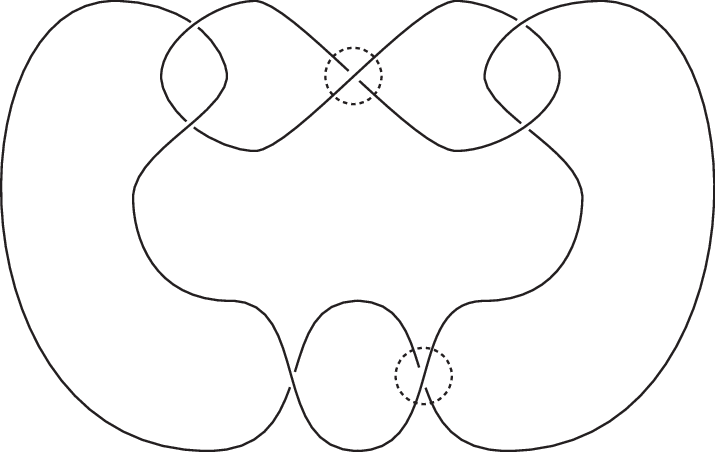}
}
\vskip 6pt
\centerline{Figure 8} 

The diagram shown in Figure 8 is rather undistinguished. The Whitehead link has an alternating
diagram, with fewer crossings, shown in Figure 9. The crossing arcs associated with the two
crossing changes are shown in Figure 9. One of these is isotopic to a vertical arc at one of the
crossings in Figure 9. The other one can be made vertical, if one first performs a flype
on the diagram, taking it to another alternating diagram. Thus, we obtain the following corollary.

\noindent {\bf Theorem 1.3.} {\sl Any crossing change that turns the Whitehead link into a split link is
equivalent to changing some crossing in some alternating diagram.}

\vskip 12pt
\centerline{
\epsfxsize = 2.3in
\epsfbox{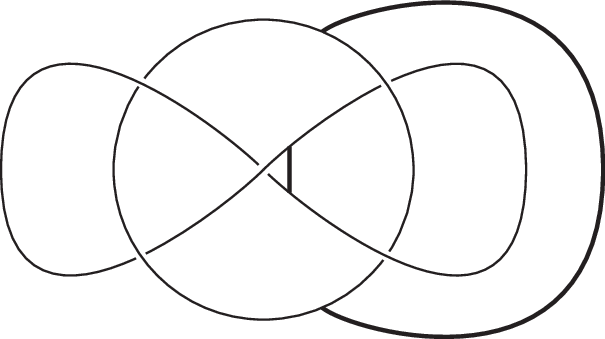}
}
\vskip 6pt
\centerline{Figure 9}

Associated to each of the two crossing changes in Figure 8, there are two crossing circles.
Two of these are isotopic to each other. Thus, we get at most three inequivalent
crossing circles in total. In fact, these are readily seen to be inequivalent, for example,
by examining their linking number with the components of the Whitehead link.

\noindent {\sl Proof of Theorem 12.1.} We follow the procedure given in Section 10. 
Note first that the Whitehead link is hyperbolic and 2-string prime. The latter fact
can be proved by observing that it is a 2-bridge link and hence its double branched cover
is a lens space, and then using Theorem 4.1.

We consider all sublinks $L'$ consisting of just one component.
Since there is an ambient isotopy that swaps the two components of the Whitehead link $L$, we may fix $L'$ to be one specific
component. The double cover of $S^3 - {\rm int}(N(L - L'))$ branched over $L'$ is shown in Figure 10. 
It is the exterior $M$ of the $(4,2)$ torus link. 

\vskip 12pt
\centerline{
\epsfxsize = 3.9in
\epsfbox{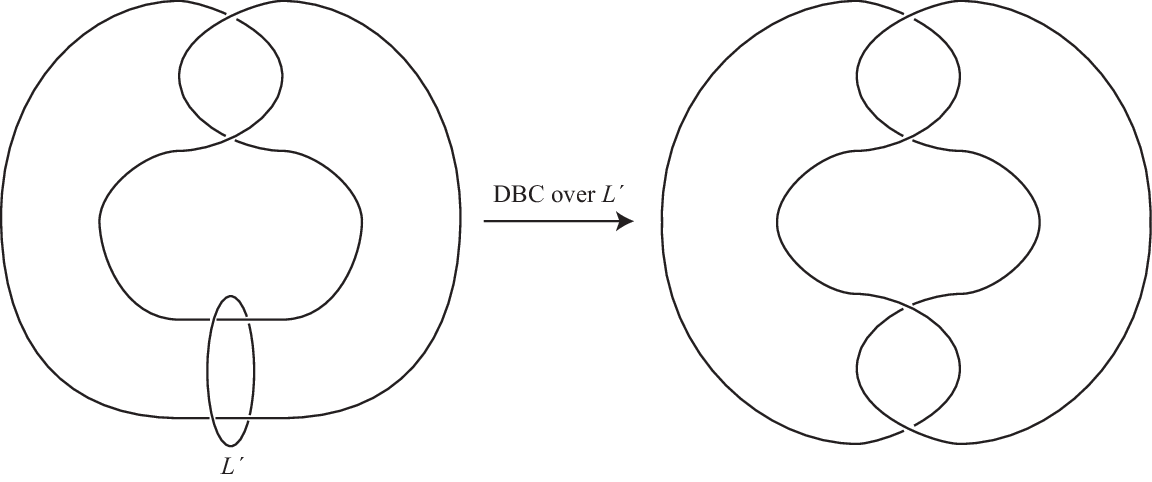}
}
\vskip 6pt
\centerline{Figure 10}

We wish to produce a finite list of knots $K$ in $M$ with slopes $\sigma$ that satisfy
the hypotheses of Theorem 5.1. Note that the arguments in Section 7 apply
and so we may assume that $M - {\rm int}(N(K))$ is irreducible and atoroidal. 

Since $M$ is Seifert fibred, we could use Theorem 6.1.
Instead, we consider the method discussed in Section 10. The relevant knots $K$
and slopes $\sigma$ arise via taut sutured manifold decompositions
$$(M, \emptyset) \buildrel S_1 \over \longrightarrow \dots \buildrel S_{n-1} \over \longrightarrow (M_{n}, \gamma_{n}),$$
where $(M_{n}, \gamma_{n})$ is a solid torus regular neighbourhood of $K$ with sutures of slope $\sigma$,
plus possibly some taut 3-balls. We may take the homology class of the surface $S_1$ to be any non-trivial class in
$H_2(M, \partial M)$ that has zero intersection number with $K$. This is because such classes in $H_2(M,\partial M)$
are precisely those in the image of the non-trivial classes in $H_2(M - {\rm int}(N(K)), \partial M)$.
We will show that we can in fact take $S_1$ to be the annulus $A$ shown in Figure 11.

\vskip 12pt
\centerline{
\epsfxsize = 2.3in
\epsfbox{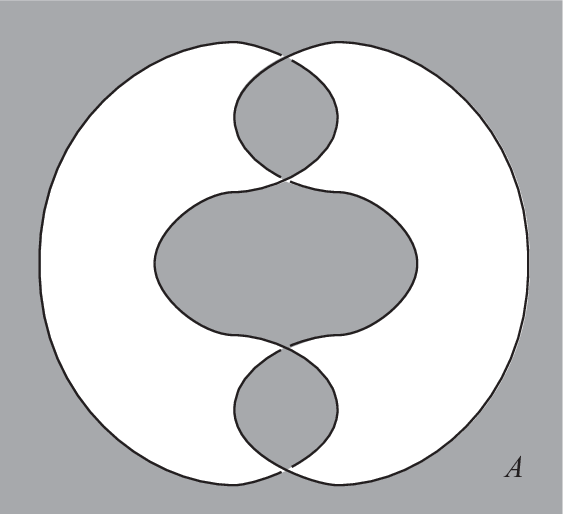}
}
\vskip 6pt
\centerline{Figure 11} 

Note first that, by the Theorem 10.1, there is an isotopy of $K$ taking it off the annulus $A$.
Hence, $K$ does have zero algebraic intersection number with $A$. So we may
take $S_1$ to be homologous to $A$. Furthermore, $S_1$ is incompressible, since it is
the first surface in a taut sutured manifold hierarchy.
We will in fact show that any connected orientable incompressible surface properly embedded in $M$
that is homologous to $A$ is isotopic to $A$. 

The manifold $M$ is Seifert fibred with base space an annulus and with a single exceptional fibre
of order $2$. Any essential properly embedded surface in $M$ is isotopic to one that is 
horizontal (that is, it is transverse to the fibres) or vertical (that is, it is a union of fibres). The annulus $A$ is
vertical. Horizontal and vertical surfaces are not homologous, as can be seen
for example by considering their algebraic intersection number with a regular fibre.
Thus, any incompressible surface homologous to $A$ is also vertical.
But, because the base space of $M$ is an annulus with a single exceptional fibre,
the unique connected orientable vertical surface that is homologically non-trivial is 
isotopic to $A$.

Thus, we may assume that $S_1$ is $A$. The second manifold $(M_2, \gamma_2)$ is therefore a  solid torus. Its boundary $T$ lies in 
$M - {\rm int}(N(K))$. Since $M - {\rm int}(N(K))$ is irreducible and atoroidal, $T$ is boundary parallel in $M - {\rm int}(N(K))$.
It is not parallel to $\partial M$, and hence it is parallel to $\partial N(K)$. Hence,
$M_2 - {\rm int}(N(K))$ is a copy of $T^2 \times I$. Therefore, $H_2(M_2 - {\rm int}(N(K)), \partial M_2)$
is trivial, and therefore $M_2$ is the final manifold in the hierarchy.
Hence, the only possibility, up to ambient isotopy, for $K$ is as shown in Figure 12.

\vskip 12pt
\centerline{
\epsfxsize = 2.2in
\epsfbox{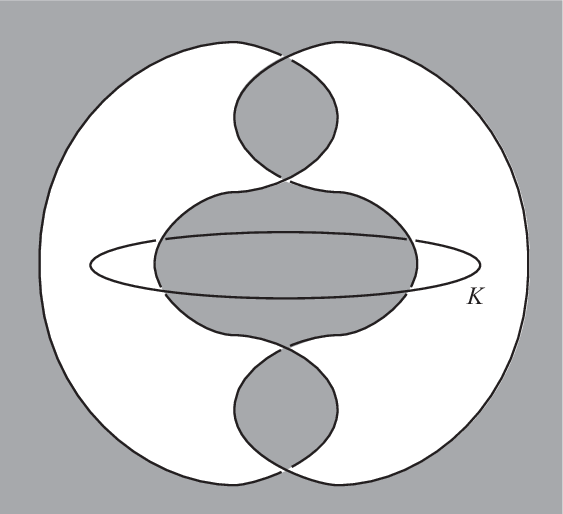}
}
\vskip 6pt
\centerline{Figure 12}

Note that $K$, as shown in Figure 12 is invariant under the involution of $M$.
The quotient arc $\beta$ is shown in Figure 13.

 \vskip 12pt
\centerline{
\epsfxsize = 1.8in
\epsfbox{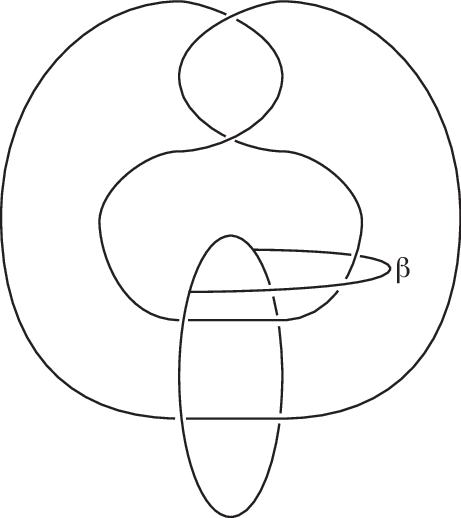}
}
\vskip 6pt
\centerline{Figure 13}

The final subroutine of the algorithm in Section 9 now produces two possibilities for the arc $\alpha$.
One of these is $\beta$. To produce the other one, we note that the
slope of $\gamma_n$ is $1/2$. Applying Theorem 2.1, we see that conclusion (ii)
there again gives $\alpha = \beta$.

Finally observe that tangle replacement along $\beta$ does change $L$
into a split link. Thus, we deduce that this is the only possibility for the arc $\alpha$ with endpoints
in $L'$. An isotopy takes $L \cup \beta$ to the link and one of the crossing arcs shown in Figure 9.
$\square$

Theorem 1.3 provides some evidence for the following conjecture.

\noindent {\bf Conjecture 12.2.} {\sl Any crossing change that turns an alternating link into a split link
is equivalent to changing some crossing in some alternating diagram.}

It is a theorem of McCoy [31] that an alternating knot has unknotting number one if and only
if one can change a crossing in some alternating diagram of the knot and obtain the unknot. 
However, this does not imply, of course, that every crossing change that turns an alternating knot
into the unknot is equivalent to changing some crossing in some alternating diagram.
Indeed it seems unlikely that the methods developed by McCoy, which use Heegaard Floer
homology, would lead to a proof of Conjecture 12.2. But McCoy's theorem does lend weight to
the conjecture.

\vskip 18pt
\centerline{\caps References}
\vskip 6pt



\item{1.} {\caps K. Baker, D. Buck,} {\sl The classification of rational subtangle replacements between rational tangles,} Algebraic and Geometric Topology 13 (2013)1413--1463.


\item{2.} {\caps R. Bing}, {\sl A homeomorphism between the {$3$}-sphere and the sum of two
              solid horned spheres}, {Ann. of Math. (2)} {56} (1952) 354--362.

\item{3.} {\caps M. Boileau, B. Leeb, J. Porti,} {\sl Geometrization of 3-dimensional orbifolds.} Ann. of Math. (2) 162 (2005), no. 1, 195--290. 

\item{4.} {\caps F. Bonahon, L. Siebenmann}, {\sl 
The characteristic toric splitting of irreducible compact 3-orbifolds},
Math. Ann. 278 (1987) 441--479.

\item{5.} {\caps F. Bonahon, L. Siebenmann}, {\sl 
New Geometric Splittings of Classical Knots and the Classification and Symmetries of Arborescent Knots,} Preprint.

\item{6.} {\caps J. Cha, S. Friedl, M. Powell,} {\sl Splitting numbers of links}, Proc. Edinb. Math. Soc. (2) 60 (2017), no. 3, 587--614.

\item{7.} {\caps D. Cooper, C. Hodgson, S. Kerckhoff,} {\sl Three-dimensional orbifolds and cone-manifolds.} With a postface by Sadayoshi Kojima. 
MSJ Memoirs, 5. Mathematical Society of Japan, Tokyo, 2000.

\item{8.} {\caps M. Culler, C. McA. Gordon, J. Luecke, P. Shalen,} {\sl Dehn surgery on knots.} Ann. of Math. (2) 125 (1987), no. 2, 237--300.

\item{9.} {\caps D. Gabai}, {\sl Foliations and the topology of 3-manifolds. III.} J. Differential Geom. 26 (1987), no. 3, 479--536.

\item{10.} {\caps D. Gabai}, {\sl Surgery on knots in solid tori.} Topology 28 (1989), no. 1, 1--6. 

\item{11.} {\caps D. Gabai}, {\sl The Smale conjecture for hyperbolic 3-manifolds: ${\rm Isom}(M^3)\simeq {\rm Diff}(M^3)$.}
 J. Differential Geom. 58 (2001), no. 1, 113--149.

\item{12.} {\caps C. Gordon, J. Luecke,} {\sl Reducible manifolds and Dehn surgery.}
Topology 35 (1996), no. 2, 385--409.

\item{13.} {\caps D. Yu. Grigor\'ev, N. Vorobjov, Jr,} {\sl Solving systems of polynomial inequalities in subexponential time,} J. Symbolic Comput. 5 (1988) 37--64. 

\item{14.} {\caps W. Haken,} {\sl Theorie der Normalfl\"achen.} Acta Math. 105 (1961) 245--375.



\item{15.} {\caps W. H. Holzmann,} {\sl An equivariant torus theorem for involutions}, Trans. Amer.
Math. Soc 326 (1991) 887--906.

\item {16.} {\caps W. Jaco,} {\sl Lectures on three-manifold topology.} 
CBMS Regional Conference Series in Mathematics, 43. American Mathematical Society, Providence, R.I., 1980.

\item{17.} {\caps W. Jaco, U. Oertel,}
{\sl An algorithm to decide if a 3-manifold is a Haken manifold.}
Topology 23 (1984), no. 2, 195--209.

\item{18.} {\caps W. Jaco, J. H. Rubinstein,} {\sl PL minimal surfaces in 3-manifolds.} J. Differential Geom. 27 (1988), no. 3, 493--524.

\item{19.} {\caps W. Jaco, J. Tollefson,} {\sl Algorithms for the complete decomposition of a closed 3-manifold},
Illinois J. Math. 39 (1995) 358--406.

\item{20.} {\caps D. A. Ka\u{z}dan, G. A. Margulis,} {\sl A proof of SelbergÕs hypothesis,} Mat. Sb. (N.S.)
75 (117) (1968), 163--168. 

\item{21.} {\caps T. Kobayashi,} {\sl Equivariant annulus theorem for 3-manifolds.} Proc. Japan Acad. Ser. A Math. Sci. 59 (1983), no. 8, 403--406.

\item{22.} {\caps  G. Kuperberg}, {\sl  Algorithmic homeomorphism of 3-manifolds as a corollary of geometrization},
Pacific J. Math. 301 (2019) 189--241. 

\item{23.} {\caps M. Lackenby,} {\sl Surfaces, surgery and unknotting operations,}
Math. Ann. 308 (1997) 615--632.

\item{24.} {\caps M. Lackenby,} {\sl Exceptional surgery curves in triangulated 3-manifolds,} Pacific J. Math. 210 (2003), 101--163.


\item{25.} {\caps W. B. R. Lickorish,} {\sl The unknotting number of a classical knot.} Contemp. Math. 44 (1985) 117--121.

\item{26.} {\caps C. Livingston,} {\sl KnotInfo}, http://www.indiana.edu/$\sim$knotinfo

\item{27.} {\caps C. Livingston,} {\sl Splitting numbers of links and the four-genus.} Proc. Amer. Math. Soc. 146 (2018), no. 1, 421--427.


\vfill\eject
\item{28.} {\caps J. Manning,} {\sl Algorithmic detection and description of hyperbolic structures on closed 3-manifolds with solvable word problem.} Geom. Topol. 6 (2002), 1--25.

\item{29.} {\caps S. Matveev,} {\sl Transformations of special spines, and the Zeeman conjecture.} Izv. Akad. Nauk SSSR Ser. Mat. 51 (1987), no. 5, 1104--1116, 1119; translation in Math. USSR-Izv. 31 (1988), no. 2, 423Ð434

\item{30.} {\caps S. Matveev}, {\sl Algorithmic topology and classification of 3-manifolds.}
Algorithms and Computation in Mathematics, 9. Springer, Berlin, 2007.

\item{31.} {\caps D. McCoy,} {\sl Alternating knots with unknotting number one.} Adv. Math. 305 (2017), 757--802. 

\item{32.} {\caps W. H. Meeks III,  L. Simon, S.-T. Yau,} {\sl 
Embedded Minimal Surfaces, Exotic Spheres, and Manifolds with Positive Ricci Curvature},
Ann. Math. 116, No. 3 (1982) 621--659.

\item{33.} {\caps W. H. Meeks III, S.-T. Yau,} {\sl The equivariant Dehn's lemma and the loop theorem,} Comment. Math.
Helv. 56 (1981), 225--239.

\item{34.} {\caps E. Moise,} {\sl Geometric topology in dimensions {$2$} and {$3$}}, {Graduate Texts in Mathematics, Vol. 47},
Springer-Verlag, New York-Heidelberg, 1977.	

\item{35.} {\caps J. Montesinos,} {\sl Three manifolds as 3-fold branched covers of $S^3$}, Quart. J. Math. Oxford (2), 27 (1976), 85--94.

\item{36.} {\caps J. Morgan,} {\sl The Smith conjecture.} (New York, 1979) Pure Appl. Math., 112.


\item{37.} {\caps Y. Nakanishi,} {\sl A note on unknotting number.} Math. Sem. Notes Kobe Univ. 9 (1981) 99--108.

\item{38.} {\caps B. Owens,} {\sl Unknotting information from Heegaard Floer homology.} Adv. Math. 217 (2008), no. 5,
2353--2376.

\item{39.} {\caps P. Ozsv\'ath, Z. Szab\'o}, {\sl Knots with unknotting number one and Heegaard Floer homology,} Topology
44 (2005), no. 4, 705--745.

\item{40.} {\caps G. Perelman,} {\sl The entropy formula for the Ricci flow and its geometric applications,} Preprint,
arxiv:math.DG/0211159

\item{41.} {\caps  G. Perelman,} {\sl Ricci flow with surgery on three-manifolds,} Preprint, arxiv:math.DG/0303109

\item{42.} {\caps G. Perelman,} {\sl Finite extinction time for the solutions to the Ricci flow on certain three-manifolds,}
Preprint, arxiv:math.DG/0307245

\item{43.} {\caps C. Petronio, J. Weeks,} {\sl Partially flat ideal triangulations of cusped hyperbolic 3-manifolds},
Osaka J. Math. 37 (2000), 453--466.

\item{44.} {\caps M. Scharlemann,} {\sl Sutured manifolds and generalized Thurston norms.} 
J. Differential Geom. 29 (1989), no. 3, 557--614.

\vfill\eject
\item{45.} {\caps M. Scharlemann,} {\sl Producing reducible 3-manifolds by surgery on a knot.} Topology 29 (1990), no. 4, 481--500.

\item{46.} {\caps P. Scott, H. Short,} {\sl The homeomorphism problem for closed 3-manifolds,} Algebr. Geom. Topol. 14 (2014), no. 4,
2431--2444

\item{47.} {\caps A. Stoimenow,} {\sl Polynomial values, the linking form and unknotting numbers.} Math. Res. Lett. 11
(2004), no. 5-6, 755--769.

\item{48.} {\caps T. Tanaka,} {\sl Unknotting numbers of quasipositive knots,} Topology Appl. 88 (1998), no. 3, 239--246.

\item{49.} {\caps A. Tarski,} {\sl A Decision Method for Elementary Algebra and Geometry,} RAND Corporation, Santa Monica, Calif., 1948.

\item{50.} {\caps W. Thurston}, {\sl The geometry and topology of three-manifolds},
library.msri.org/books/gt3m

\item{51.} {\caps J. Tollefson,} {\sl Involutions of Seifert fiber spaces,} Pacific J. Math. 74 (1978), no. 2, 519--529.

\item{52.} {\caps J. Weeks,} {\sl Snappea,} http://www.geometrygames.org/SnapPea/

\item{53.} {\caps J. Weeks,} {\sl Convex hulls and isometries of cusped hyperbolic 3-manifolds,} Topology and its Applications 52 (1993), no. 2, 127--149.


\end